\documentclass[12pt,reqno]{amsart}
\usepackage{amssymb}
\usepackage{cases}
\usepackage{bbm}
\usepackage{mathrsfs}
\usepackage{graphicx}
\usepackage{amsfonts}
\usepackage{enumerate}

\usepackage{amssymb, amsmath, amsfonts, latexsym}
\usepackage{enumerate}
\usepackage[notref,notcite]{showkeys}

\allowdisplaybreaks[4]

\newtheorem{thm}{Theorem}[section]
\newtheorem{corollary}[thm]{Corollary}
\newtheorem{lemma}[thm]{Lemma}
\newtheorem{proposition}[thm]{Proposition}
\newtheorem{example}[thm]{Example}
\newtheorem{remark}[thm]{Remark}

\newtheorem{hyp}[thm]{Hypothesis}

\theoremstyle{definition}
\newtheorem{definition}{Definition}[section]
 \theoremstyle{remark}

\numberwithin{equation}{section}

\def\<{\langle}
\def\>{\rangle}

\def\beq{\begin{equation}}
\def\nneq{\end{equation}}

\def\bdef{\begin{defn}}
\def\ndef{\end{defn}}

\def\bthm{\begin{thm}}
\def\nthm{\end{thm}}

\def\bprop{\begin{proposition}}
\def\nprop{\end{proposition}}

\def\brmk{\begin{remarks}}
\def\nrmk{\end{remarks}}

\def\bexa{\begin{exa}}
\def\nexa{\end{exa}}

\def\blem{\begin{lemma}}
\def\nlem{\end{lemma}}

\def\bcor{\begin{corollary}}
\def\ncor{\end{corollary}}

\def\benu{\begin{enumerate}}
\def\nenu{\end{enumerate}}

\def\bhyp{\begin{hyp}}

\def\nhyp{\end{hyp}}

\def\<{\langle}
\def\>{\rangle}

\date{}

\def\bexe{\begin{exe}}
\def\nexe{\end{exe}}

\def\bprf{\begin{proof}}
\def\nprf{\end{proof}}

\def\bdes{\begin{description}}
\def\ndes{\end{description}}

\title[Solutions of MMVSDEs driven by jump process]{Existence of Solutions for Multivalued Mckean-Vlasov SDEs with Non-Lipschitz Coefficients driven by Jump Processes}

\author{Lingyan Cheng}
\address{Lingyan Cheng. School of Mathematics and Statistics, Nanjing University of Science and Technology,
Nanjing 210094, Jiangsu, PR China.}
\email{cly@njust.edu.cn}

\author{Caihong Gu}
\address{Caihong Gu. School of Mathematics and Statistics, Wuhan University, Wuhan 430072, Hubei, PR China.}
\email{gucaihong@whu.edu.cn}

\author{Wei Liu}
\address{Wei Liu\\School of Mathematics and Statistics, Wuhan University, Wuhan 430072, Hubei, PR China.}
\email{wliu.math@whu.edu.cn}

\author{Fengwu Zhu}
\address{Fengwu Zhu. School of Mathematics and Statistics, Wuhan University, Wuhan 430072, Hubei, PR China.}
\email{fwzhu\_math@whu.edu.cn}

\date{}

\begin{document}
\maketitle

\noindent {\bf Abstract:}
In this paper, we first establish the existence and uniqueness of strong solutions for multivalued McKean-Vlasov stochastic differential equations (MMVSDEs) driven by L\'evy noise with non-Lipschitz coefficients. It is important to note that these findings are based upon the well-posedness of strong solutions for MMVSDEs under Lipschitz conditions, which will be stated briefly. Secondly, we study the existence of weak solutions under linear growth condition. Finally, we prove the existence of martingale solutions.
 \vskip0.3cm

 \vskip0.3cm

\vskip0.3cm

\section{Introduction}\label{sec1}
\setcounter{equation}{0}
\noindent

Consider the following multivalued McKean-Vlasov stochastic differential equation on $\mathbb{R}^d$:
\begin{eqnarray}\label{equ1.1}
    dX_t\in -A(X_t)dt+b[X_t,\mathcal{L}_{X_t}]dt+\sigma[X_t,\mathcal{L}_{X_t}]dW_t\nonumber\\
        \quad\quad\quad+\int_{\mathbb{Z}}G[X_{t-},\mathcal{L}_{X_t},z]\widetilde{N}(dt,dz),\quad t\in[0,1],
\end{eqnarray}
where $A$ is a multivalued maximal monotone operator on $\mathbb{R}^d$ (see Definition 2.1), $\mathcal{L}_{X_t}$ denotes the law of $X_t$, and $b$, $\sigma$, $G$ are the functionals of $X_t$ and $\mathcal{L}_{X_t}$. For a given probability measure $\mu$, denote
\begin{eqnarray*}
    b[x,\mu]:=\int_{\mathbb{R}^d} b(x,y)\mu(dy),~~\sigma [x,\mu]:=\int_{\mathbb{R}^d} \sigma(x,y)\mu(dy),\\
    \int_{\mathbb{Z}}G[x,\mu,z]\nu(dz):=\int_{\mathbb{Z}}\int_{\mathbb{R}^d} G(x,y,z)\mu(dy)\nu(dz),
\end{eqnarray*}
where $b: \mathbb{R}^d\times\mathbb{R}^d\rightarrow\mathbb{R}^d$, $\sigma: \mathbb{R}^d\times\mathbb{R}^d\rightarrow\mathbb{R}^{d}\otimes\mathbb{R}^d$ and $G: \mathbb{R}^d\times\mathbb{R}^d\times\mathbb{Z}\rightarrow\mathbb{R}^d$ are continuous functions. Here $\mathbb{Z}$ is a locally compact Polish space with Borel $\sigma$-algebra $\mathcal{B}(\mathbb{Z})$. $W$ is a standard $d$-dimensional Brownian motion and $N$ is a Poisson random measure with intensity measure $\nu$ with $\nu(\mathbb{Z}\backslash B)<+\infty$ for all $B\in\mathcal{B}(\mathbb{Z})$, on some filtered probability space $(\Omega, \mathcal{F}, \mathbb{P})$. The compensated Poisson measure is given by $\widetilde{N}(dt,dz):=N(dt,dz)-\nu(dz)dt$. $W$ and $N$ are assumed to be independent. Without loss of generality, our model is limited in the time interval $[0,1]$.

In the case $A=G=0$, Eq.(\ref{equ1.1}) is the classical McKean-Vlasov stochastic differential equation (MVSDE) driven by Brownian motion. These equations were initially proposed by  Kac \cite{Kac1, Kac2} as a stochastic toy model for the Vlasov kinetic equation of plasma, and later introduced by McKean \cite{Mckean} to model plasma dynamics. One can refer to \cite{Cattiaux, Graham, Guerin, Guillin1, Guillin2, Liu1, Liu2, Mckean, Mishura} and the reference therein. There has been an increasing interest to study the existence and uniqueness for the solutions of MVSDEs. Wang \cite{Wang} established the strong well-posedness of MVSDEs with one-sided Lipschitz continuous drifts. Under integrability conditions on distribution dependent coefficients, R\"{o}ckner et al. \cite{Huang,Rockner} showed the strong well-posedness of MVSDEs with Lipschitz continuous coefficients. Mehri and Stannat \cite{Mehri} presented a Lyapunov-type approach to prove the existence and uniqueness for the solutions of MVSDEs.

When $A\neq0$ and $G=0$, in the distribution independent case, Eq.(\ref{equ1.1}) is the multivalued stochastic differential equation (MSDE) driven by Brownian motion. C\'epa \cite{Cepa,Cepa1} first introduced a pair of continuous $\mathcal{F}_t$-adapted processes $(X_t,K_t)$ to solve the MSDE. Z\v{a}linescu \cite{Zalinescu} proved the existence of weak solution to  finite dimensional MSDE by martingale approach. Wu \cite{Wu} proved that a non-explosion solution exists when the drift coefficient satisfies linear growth conditions and the diffusion coefficient is uniformly elliptic. For more details, the reader is referred to \cite{Cepa2,Ren1,Ren2,Zhang,Hu}. On the other hand, for the distribution dependent case, Chi \cite{Chi} proved the existence and uniqueness of strong solutions under global lipschitz conditions, and obtained the existence of the weak solutions  under continuous and linear growth assumptions. Recently, Gong and Qiao \cite{Gong} established the well-posedness and stability under non-Lipschitz conditions.

When $A\neq0$ and $G\neq0$, there are few results about the existence and uniqueness for the solutions of MSDEs with jumps. In the distribution independent case, Ren and Wu \cite{Ren3} proved the existence and uniqueness for the solutions of MSDEs driven by Poisson point processes under Lipschitz continuous assumptions on the coefficients, when the domain of the multivalued operator is the whole space $\mathbb{R}^d$.

However, as far as we know in the jump case, the solutions of MSDEs with non-Lipschitz coefficients, as well as those of MMVSDEs, have not been studied. The present paper is devoted to investigating the existence and uniqueness of strong solutions for general MMVSDEs with non-Lipschitz coefficients driven by jump process, which extend the results mentioned before. The most important difference between the MVSDEs with jumps and without jumps is the tightness of the solutions for a sequence of equations. The lack of continuity for solution makes us to find a suitable modulus of continuity to characterize the tightness.

Unlike the classical stochastic differential equations, a pair of $\mathcal{F}_t$-adapted process $(X_t, K_t)$ is called a solution of Eq.(\ref{equ1.1}) if
\begin{eqnarray*}
    X_t=X_0-K_t+\int_0^t b[X_s,\mathcal{L}_{X_s}]ds+\int_0^t \sigma[X_s,\mathcal{L}_{X_s}]dW_s \\
    ~~~~+\int_0^t \int_{\mathbb{Z}} G[X_{s-},\mathcal{L}_{X_s},z]\widetilde{N}(ds,dz),~ t\in[0,1],
\end{eqnarray*}
where $K_t$ is a locally finite variation process and satisfies certain conditions. When the coefficients $b$, $\sigma$ and $G$ are all globally Lipschitz continuous, the Picard iteration can be employed to prove the existence and uniqueness of a strong solution to Eq.(\ref{equ1.1}). When the coefficients satisfy non-Lipschitz continuity and linear growth condition, the problem may be more involved, as discussed in \cite{Cattiaux}. In this  paper we overcome this difficulty by an application of smoothing technique, which allows us to find a sequence of Lipschitz functions to approximate the non-Lipschitz coefficients, through constructing convolutions of the non-Lipschitz functions and suitable kernels.

In fact, the equivalence between the existence of weak solution and that of martingale solution in the distribution dependent case can be established by using fixed distribution methods and the results in Kurtz \cite{Kurtz}. It is important to note that martingale solution and weak solution are essentially two different viewpoints of the same underlying concept. The former emphasizes the probabilistic structure (via generator/martingale conditions), while the latter is more concerned with construction of stochastic processes and the modeling of noise, highlighting a different but equally essential aspect of the solution framework. Hence we prefer to prove the existence of martingale solution and weak solution for Eq.(\ref{equ1.1}) separately.

The rest of this paper is organized as follows. In Section $2$, we introduce the necessary notations and present our main results. Then we establish the existence and uniqueness of the strong solutions for Eq.$(\ref{equ1.1})$ under globally Lipschitz conditions and then non-Lipschitz conditions in Section $3$. Section $4$ is devoted to investigating the existence of weak solutions under the non-Lipschitz assumptions.  The existence of martingale solution is proved by use of the tightness of probability measures in the last Section.

\section{Preliminaries and Main Results}\label{sec2}
\setcounter{equation}{0}

In this section,we will introduce the notations and recall some known results which will be used later in the proofs.
\subsection{Multivalued Maximal Monotone Operators}

Throughout this paper, we use the following notations and definitions:

$\bullet$ $\mathcal{B}(S)$: The Borel $\sigma$-algebra on a metric space $S$.

$\bullet$ $\mathbb{C}([0,1],S)$:  The space of continuous functions from the interval $[0,1]$ to the metric space $S$, equipped with the uniform convergence topology.

$\bullet$ $\mathbb{D}([0,1],S)$: The space of $c\grave{a}dl\grave{a}g$ (right continuous with left limit) functions from the interval $[0,1]$ to the metric space $S$, equipped with the Skorohod topology.

For the sake of simplicity, we denote by $|\cdot|$, $\|\cdot\|$, $\|\cdot\|_{\mathbb{C}}$ and $\|\cdot\|_{\mathbb{D}}$, the Euclidean norm in $\mathbb{R}^d$, the Hilbert-Schmidt norm of matrix in $\mathbb{R}^d\otimes\mathbb{R}^d$, the suprmum norm in $\mathbb{C}:=\mathbb{C}([0,1],\mathbb{R}^d)$ and the supremum norm in $\mathbb{D}:=\mathbb{D}([0,1],\mathbb{R}^d)$, respectively. For a Polish space $\mathbb{X}$, denote $\mathcal{P}(\mathbb{X})$ as the space of all probability measures on Borel $\sigma$-algebra $\mathcal{B}(\mathbb{X})$ carrying the usual topology of weak convergence. $\mathcal{P}_p(\mathbb{X})$ denotes the subset of $\mathcal{P}(\mathbb{X})$ with finite p-order moment.

Denote $2^{\mathbb{R}^d}$ the collection of all subsets of $\mathbb{R}^d$. A $multivalued$ $operator$ $A$: $\mathbb{R}^d\rightarrow 2^{\mathbb{R}^d}$ can be identified by its $graph$
\begin{eqnarray*}
    \mathcal{G}(A):=\Big\{(x,x^*)\in\mathbb{R}^d\times\mathbb{R}^d|x\in \mathcal{D}(A), x^*\in A(x)\Big\},
\end{eqnarray*}
where $\mathcal{D}(A):=\{x\in\mathbb{R}^d\mid A(x)\neq\emptyset\}$ and $\text{Im}(A):=\cup_{x\in \mathcal{D}(A)} A(x)$ denote the domain and image of $A$, respectively. We refer the reader to \cite{Brezis} and related literature for more details.

\begin{definition}
A multivalued operator $A$ is called monotone if
\begin{eqnarray*}
    \langle x-x',y-y'\rangle\geq0,~~\forall (x,y),(x',y')\in \mathcal{G}(A).
\end{eqnarray*}
A monotone operator $A$ is called maximal if
\begin{eqnarray*}
(x,y)\in \mathcal{G}(A)\Leftrightarrow\langle x-x',y-y'\rangle\geq0,~~\forall(x',y')\in \mathcal{G}(A),
\end{eqnarray*}
which is equivalent to (cf. Proposition 2.2 in \cite{Brezis})
\begin{eqnarray*}
    \text{Im}(I+A)=\mathbb{R}^d.
\end{eqnarray*}
\end{definition}

If $A$ is a multivalued maximal monotone operator on $\mathbb{R}^d$, then $\text{Int} (\mathcal{D}(A))$ and $\overline{\mathcal{D}(A)}$ are both convex sets, and
\begin{eqnarray*}
    \text{Int}(\mathcal{D}(A))=\text{Int}(\overline{\mathcal{D}(A)}).
\end{eqnarray*}
Moreover, for every $x\in\mathbb{R}^d$, $A(x)$ is a closed convex set of $\mathbb{R}^d$. A crucial property of a maximal monotone operator $A$ is its local boundedness on $\text{Int}(\mathcal{D}(A))$, which means that the union $\cup_{x\in\mathbb{K}}A(x)$ is bounded for every compact set $\mathbb{K}\subseteq \text{Int}(\mathcal{D}(A))$.

Next we give two examples to illustrate the multivalued maximal monotone operators.
\begin{example}
For a lower semicontinuous convex function $\varphi:\mathbb{R}^d\mapsto(-\infty,+\infty]$, let $\mathcal{D}(\varphi)\equiv\{x\in\mathbb{R}^d:\varphi(x)<+\infty\}$ and $\text{Int}(\mathcal{D}(\varphi))$ be the interior of $\mathcal{D}(\varphi)$. Assume $\text{Int}(\mathcal{D}(\varphi))\neq\emptyset$. Define the sub-differential operator of $\varphi$:
\begin{eqnarray*}
    \partial\varphi(x):=\Big\{x^*\in\mathbb{R}^d\mid\langle y-x,x^*\rangle+\varphi(x)\leq\varphi(y),\forall y\in\mathbb{R}^d\Big\}.
\end{eqnarray*}
Then it is easy to verify that $\partial\varphi$ is a maximal monotone operator.
\end{example}

\begin{example}
For a closed convex subset $\mathcal{O}$ of $\mathbb{R}^d$, suppose $\text{Int}(\mathcal{O})\neq\emptyset$. Define the indicator function of $\mathcal{O}$ as follows:
\begin{eqnarray*}
    \text{I}_{\mathcal{O}}(x):=\Big\{\begin{array}{ll}
                            0,\ \ \ \ \ \ \ if\ x\in\mathcal{O}, \\
                            +\infty,\ \ \ if\ x\in\mathbb{R}^d\backslash\mathcal{O}.
                          \end{array}
\end{eqnarray*}
The sub-differential operator of $\text{I}_{\mathcal{O}}$ is given by
\begin{eqnarray*}
    \partial \text{I}_{\mathcal{O}}(x):=\Bigg\{\begin{array}{lll}
                                 {0},\ \ \ \ \ \ \ \ if\ x\in \text{Int}(\mathcal{O}), \\
                                 \Pi_x,\ \ \ \ \ \ if\ x\in \partial(\mathcal{O}), \\
                                 \emptyset,\ \ \ \ \ \ \ \ if\ x\in\mathbb{R}^d\backslash\mathcal{O},
                               \end{array}
\end{eqnarray*}
where $\Pi_x$ is the exterior normal cone at $x$. It is also easy to verify that $\text{I}_{\mathcal{O}}$ is a maximal monotone operator.
\end{example}


Define $\mathcal{V}_0$ as the set of all continuous functions $K: [0,1]\rightarrow\mathbb{R}^d$ with finite variation and $K_0=0$. For any $K\in\mathcal{V}_0$ and $t\in[0,1]$, denote the variation of $K$ on $[0,t]$ by $\|K\|_0^t$ and $\|K\|_{TV}:=\|K\|_0^1$. Set
\begin{eqnarray*}
    \mathcal{A}&:=&\Big\{(X,K)\mid X\in\mathbb{D}([0,1],\overline{\mathcal{D}(A)}),K\in\mathcal{V}_0,\forall(x,y)\in \mathcal{G}(A),\\
    &&\quad\quad\quad\quad\langle X_t-x,dK_t-ydt\rangle\geq0\Big\}.
\end{eqnarray*}

The proof of the following proposition is analogous to that presented in C\'epa \cite{Cepa1}, so we omit it here for brevity.

\begin{proposition}\label{pro2.1}
Let $(X,K)$ be a pair of functions such that $X\in\mathbb{D}([0,1],\overline{\mathcal{D}(A)})$ and $K\in\mathcal{V}_0$. The following statements are equivalent:
\begin{enumerate}[(i)]
	\item $(X,K)\in\mathcal{A}$.

\item for any $(x,y)\in\mathbb{D}([0,1],\mathbb{R}^d)\times \mathbb{C}([0,1],\mathbb{R}^d)$ with $(x_t,y_t)\in \mathcal{G}(A)$, it holds that
\begin{eqnarray*}
    \langle X_t-x_t,dK_t-y_tdt\rangle\geq0.
\end{eqnarray*}

\item for any $(X',K')\in\mathcal{A}$, it holds that
\begin{eqnarray*}
    \langle X_t-X'_t,dK_t-dK'_t\rangle\geq0.
\end{eqnarray*}
\end{enumerate}
\end{proposition}

We recall the following results, which are essentially due to C\'epa \cite{Cepa1} (see also \cite{Zhang}).

\begin{lemma}\label{lem2.1}
Assume that $\text{Int}(\mathcal{D}(A))$ is non-empty. Given any $a_0\in \text{Int}(\mathcal{D}(A))$, there exist constants $\kappa_1>0$, $\kappa_2\geq0$ and $\kappa_3\geq0$ such that for any pair of functions $(X,K)\in\mathcal{A}$ and $0\leq s<t\leq 1$, the following inequality holds:
\begin{eqnarray*}
    \int_s^t\langle X_r-a_0,dK_r\rangle\geq\kappa_1\|K\|_s^t-\kappa_2\int_s^t|X_r-a_0|dr-\kappa_3(t-s).
\end{eqnarray*}
\end{lemma}

\begin{lemma}\label{lem2.2}
For a sequence of functions $\{K^n,n\in\mathbb{N}\}\subset\mathcal{V}_0$, if $\sup_{n\in\mathbb{N}}\|K^n\|_{TV}<+\infty$ and $K^n$ converges to some $K$ in $\mathbb{C}([0,1],\mathbb{R}^d)$ as $n\rightarrow+\infty$, then $K\in\mathcal{V}_0$ and for any $X^n\in\mathbb{D}([0,1],\mathbb{R}^d)$ that converges to $X$ in $\mathbb{D}([0,1],\mathbb{R}^d)$,
\begin{eqnarray*}
   \lim_{n\rightarrow+\infty}\int_0^1\langle X_s^n,dK_s^n\rangle=\int_0^1\langle X_s,dK_s\rangle.
\end{eqnarray*}
\end{lemma}


\subsection{Solutions to Multivalued MVSDEs with Jumps}

Let $(\Omega,\mathcal{F},\mathbb{P})$ be a complete probability space and $W$ a $d$-dimensional Brownian motion defined on it. Fix $\mathbb{Z}_0=\{|z|\leq\alpha\}$ where $\alpha>0$ is a constant, with $\nu(\mathbb{Z}\backslash\mathbb{Z}_0)<+\infty$ and $\int_{\mathbb{Z}_0}|z|^2\nu(dz)<+\infty$. Let $\zeta_t$ be a stationary Poisson point process with characteristic measure $\nu$ and $D_{\zeta_t}$ be a countable subset of $[0,+\infty)$. Denote by $N(dt,dz)$ the Poisson random measure associated with $\zeta_t$, i.e.
\begin{eqnarray*}
    N(t,A):=\sum_{s\in D_{\zeta_s},s\leq t}I_A(\zeta_s)
\end{eqnarray*}
with intensity $\mathbb{E}N(dt,dz)=\nu(dz)dt$. Define $\widetilde{N}(dt,dz):=N(dt,dz)-\nu(dz)dt$ as the compensated martingale measure. Moreover, $W$ and $N$ are assumed to be independent.

Let $\{\mathcal{F}_t\}_{t\in[0,1]}$ be the filtration generated by $(W_t)_{t\geq0}$ and $N(t,\cdot)$, augmented by a $\sigma$-field $\mathcal{F}^0$, i.e.
\begin{eqnarray*}
    \mathcal{F}_t:=(\cap_{s>t}\mathcal{F}_s^0)\vee\mathcal{F}^0,~t\in[0,1],
\end{eqnarray*}
\begin{eqnarray*}
    \mathcal{F}_t^0:=\sigma\Big\{W_s,N(s,A):0\leq s\leq t, A\in\mathcal{B}(\mathbb{Z})\Big\},
\end{eqnarray*}
where $\mathcal{F}^0\subset\mathcal{F}$ is the sets of all $\mathbb{P}$-null set. $(\Omega,\mathcal{F},\{\mathcal{F}_t\}_{t\in[0,1]},\mathbb{P})$ forms a stochastic basis.

Assume that the functionals $b:\mathbb{D}\times\mathcal{P}(\mathbb{D})\rightarrow\mathbb{R}^d$, $\sigma:\mathbb{D}\times\mathcal{P}(\mathbb{D})\rightarrow\mathbb{R}^d$ and $G:\mathbb{D}\times\mathcal{P}(\mathbb{D})\times\mathbb{Z}\rightarrow\mathbb{R}^d$
are progressively measurable with respect to the filtration $\{\mathcal{F}_t\}_{t\in[0,1]}$.

\begin{definition}\label{def2.2}
$(\mathbf{Strong~Solution})$ A pair of $\mathcal{F}_t$-adapted processes $(X_.,K_.)$ is called a strong solution of Eq.(\ref{equ1.1}) with the initial value $\xi$, if $(X_.,K_.)$ is defined on a filtered probability space $(\Omega,\mathcal{F},\{\mathcal{F}_t\}_{t\in[0,1]},\mathbb{P})$ such that

\begin{enumerate}[(i)]
\item $\mathbb{P}(X_0=\xi)=1$.

\item $\Big(X_.(\omega),K_.(\omega)\Big)\in\mathcal{A}$, $a.s.$.

\item it holds that
\begin{eqnarray*}
\mathbb{P}\Big\{\int_0^1\Big(|b[X_t,\mathcal{L}_{X_t}]|+\|\sigma[X_t,\mathcal{L}_{X_t}]\|^2+\int_{\mathbb{Z}}|G[X_{t-},\mathcal{L}_{X_t},z]|^2\nu(dz)\Big)dt<+\infty\Big\}=1
\end{eqnarray*}
and
\begin{eqnarray*}
    X_t=X_0-K_t+\int_0^t b[X_s,\mathcal{L}_{X_s}]ds+\int_0^t \sigma[X_s,\mathcal{L}_{X_s}]dW_s \\
    ~~~~+\int_0^t \int_{\mathbb{Z}} G[X_{s-},\mathcal{L}_{X_s},z]\widetilde{N}(ds,dz),~ t\in[0,1],~a.s..
\end{eqnarray*}
\end{enumerate}
\end{definition}

\begin{definition}
$(\mathbf{Weak~Solution})$ We say that Eq.(\ref{equ1.1}) admits a weak solution with initial law $\mathcal{L}_{X_0}\in\mathcal{P}(\mathbb{R}^d)$ if there exists a stochastic basis $\mathcal{S}:=(\Omega,\mathcal{F},\{\mathcal{F}_t\}_{t\in[0,1]},\mathbb{P})$, a $d$-dimensional standard $\mathcal{F}_t$-Brownian motion $(W_t)_{t\geq0}$, a compensated Poisson measure $\widetilde{N}$ and a pair of $\mathcal{F}_t$-adapted processes $(X,K)$ defined on $\mathcal{S}$ such that
\begin{enumerate}[(i)]
\item $X_0$ has the law $\mathcal{L}_{X_0}$, and $(X_.(\omega),K_.(\omega))\in\mathcal{A}$, $a.s.$.

\item it holds that
\begin{eqnarray*}
\mathbb{P}\Big\{\int_0^1\Big(|b[X_t,\mathcal{L}_{X_t}]|+\|\sigma[X_t,\mathcal{L}_{X_t}]\|^2+\int_{\mathbb{Z}}|G[X_{t-},\mathcal{L}_{X_t},z]|^2\nu(dz)\Big)dt<+\infty\Big\}=1
\end{eqnarray*}
and
\begin{eqnarray*}
    X_t=X_0-K_t+\int_0^t b[X_s,\mathcal{L}_{X_s}]ds+\int_0^t \sigma[X_s,\mathcal{L}_{X_s}]dW_s \\
    ~~~~+\int_0^t \int_{\mathbb{Z}} G[X_{s-},\mathcal{L}_{X_s},z]\widetilde{N}(ds,dz),~ t\in[0,1],~a.s..
\end{eqnarray*}
\end{enumerate}
\end{definition}

Such a solution will be denoted by $(\mathcal{S};W,\widetilde{N},X,K)$.

\begin{definition}
$(\mathbf{Uniqueness\ in\ Law})$ Let $(\mathcal{S};W,\widetilde{N},X,K)$ and $(\mathcal{S}';W',\widetilde{N}',X',K')$ be two weak solutions with $\mathcal{L}_{X_0}=\mathcal{L}_{X'_0}$. If $(X,K)$ and $(X',K')$ have the same law in the space $(\mathbb{D}\times\mathbb{C},\mathcal{B}(\mathbb{D}\times\mathbb{C}))$, then we say that the uniqueness in law holds for Eq.(\ref{equ1.1}).
\end{definition}

\begin{definition}
$(\mathbf{Pathwise~Uniqueness})$ Let $(\mathcal{S};W,\widetilde{N},X,K)$ and $(\mathcal{S}';W',\widetilde{N}',X',K')$ be two weak solutions with the same initial value $X_0=X'_0$. If $(X_t,K_t)=(X'_t,K'_t)$ for all $t\in[0,1]$, then we say that the pathwise uniqueness holds for Eq.(\ref{equ1.1}).
\end{definition}

Introduce the following notations:
\begin{eqnarray*}
&&\bar{\Omega}:=\mathbb{D}\times\mathbb{C},\ \ \bar{\mathcal{F}}:=\mathcal{B}(\mathbb{D})\times\mathcal{B}({\mathbb{C}})=\mathcal{B}(\bar{\Omega}),\\ &&\bar{\mathcal{F}}_t:=\mathcal{B}(\mathbb{D}([0,t],\mathbb{R}^d))\times\mathcal{B}(\mathbb{C}([0,t],\mathbb{R}^d)),\forall t\in[0,1],\\ &&\bar{\Omega}_0:=\{(x,\eta)\in\bar{\Omega}|\eta\in\mathcal{V}_0,\langle x(r)-z,d\eta(r)-z^*dr\rangle\geq0,~\forall(z,z^*)\in \mathcal{G}(A)\}.
\end{eqnarray*}

 Let $R$ be a probability measure defined on the measurable space $(\bar{\Omega},\bar{\mathcal{F}})$ such that $R(\bar{\Omega}_0)=1$. Denote by  $\mathcal{N}_R$ the family of $R$-negligible subsets in $\bar{\mathcal{F}}$ and $\bar{\mathcal{F}}^R:=\bar{\mathcal{F}}\vee\mathcal{N}_R$, $\bar{\mathcal{F}}^R_t:=\cap_{s>t}(\bar{\mathcal{F}}_s\vee\mathcal{N}_R)$ for $t\geq0$. Obviously $(\bar{\Omega},\bar{\mathcal{F}}^R,R,\{\bar{\mathcal{F}}^R_t\}_{t\geq0})$ is a stochastic basis with $\bar{\Omega}_0\in\bar{\mathcal{F}}^R_0$.
Hence the canonical processes $\bar{X}:\bar{\Omega}\rightarrow\mathbb{D}$, $\bar{K}:\bar{\Omega}\rightarrow\mathbb{C}$, defined by
\begin{eqnarray}\label{equ5.8}
    \bar{X}(x,\eta):=x,~~\bar{K}(x,\eta):=\Big\{\begin{array}{ll}
                                                  \eta,~~~~(x,\eta)\in\bar{\Omega}_0,\\
                                                  0,~~~~(x,\eta)\in\bar{\Omega}\setminus \bar{\Omega}_0
                                                \end{array}
\end{eqnarray}
are progressively measurable with respect to this basis.

Consider the differential operator defined by
\begin{eqnarray*}
    \mathfrak{L}(\mu)f(x):=\frac{1}{2}tr\Big(\sigma\sigma^*[x,\mu]\nabla^2\Big)f(x)+\langle b[x,\mu],\nabla f(x)\rangle,\ \ (x,\mu)\in\mathbb{D}\times\mathcal{P}(\mathbb{D})
\end{eqnarray*}
for any $f\in C^2(\mathbb{R}^d)$. Here, ${\sigma}^*$ denote the transpose of $\sigma$ and $\nabla:=\Big(\frac{\partial}{\partial x_1},\cdots,\frac{\partial}{\partial x_d}\Big)$, $\nabla^2:=\Big(\frac{\partial^2}{\partial x_i\partial x_j}\Big)_{1\leq i,j\leq d}$.

\begin{definition}
$(\mathbf{Martingale~Solution})$ A probability measure $R$ on $\bar{\Omega}$ is called a martingale solution of Eq.(\ref{equ1.1}) with initial law $\mathcal{L}_{X_0}\in\mathcal{P}_2(\overline{\mathcal{D}(A)})$, if $R(\bar{\Omega}_0)=1$ and for each $f\in C_c^2(\mathbb{R}^d)$,
\begin{eqnarray*}
    \bar{M}_t^f&:=&f(\bar{X}_t)-f(\bar{X}_0)+\int_0^t \langle \nabla f(\bar{X}_s), d\bar{K}_s\rangle-\int_0^t \mathfrak{L}(\mathcal{L}_{\bar{X}_t}) f(\bar{X}_s)ds\\
    &&-\frac{1}{2}\int_0^t\int_{\mathbb{Z}}G^*[\bar{X}_{s-},\mathcal{L}_{\bar{X}_s},z]\nabla^2f(\bar{X}_{s-})G[\bar{X}_{s-},\mathcal{L}_{\bar{X}_s},z]\nu(dz)ds,~t\in[0,1]
\end{eqnarray*}
is a $\bar{\mathcal{F}}^R$-martingale under probability space $(\bar{\Omega},\bar{\mathcal{F}}^R,R,\{\bar{\mathcal{F}}^R_t\}_{t\geq0})$, where $\mathcal{L}_{\bar{X}_s}:=R\circ \bar{X}_s^{-1}$ denotes the law of $\bar{X}_s$ under $R$.
\end{definition}

\subsection{Hypotheses}

For any $\mu, \nu \in \mathcal{P}_p(\mathbb{X})$ ($p\geq1$), define the p-order Wasserstein distance between $\mu$ and $\nu$ as follows:
\begin{eqnarray*}
    \mathcal{W}_p(\mu,\nu):=\inf_{\pi\in\Pi(\mu,\nu)}\Big(\int_{\mathbb{R}^d\times\mathbb{R}^d}|x-y|^p\pi(dx,dy)\Big)^{1/p},
\end{eqnarray*}
where $\Pi(\mu,\nu)$ denotes the set of all probability measure on $\mathbb{R}^d\times\mathbb{R}^d$ with marginal distributions $\mu$ and $\nu$. Equipped with the topology induced by Wasserstein distance, $\mathcal{P}_p(\mathbb{X})$ is a Polish space (cf. \cite{Villani}).


\begin{remark}\label{rmk2.1}
From the above definition, it is easy to deduce that for any $\mathbb{R}^d$-valued random variables $\xi$, $\eta$,
\begin{eqnarray*}
    \mathcal{W}_2(\mathcal{L}_{\xi},\mathcal{L}_{\eta})\leq\Big(\mathbb{E}|\xi-\eta|^2\Big)^{1/2},
\end{eqnarray*}
where $L_{\xi}$ and $L_{\eta}$ are the laws of $\xi$ and $\eta$ respectively.
\end{remark}

Now we give our assumptions.

\begin{hyp}

\noindent{\bf (H1)} $\text{Int}(\mathcal{D}(A))\neq\emptyset$.

\noindent{\bf (H2)} the functions $b$, $\sigma$ and $G$ are continuous in $(x,y)$, and there exists some constant $L_1>0$ such that for all $(x,y)\in\mathbb{R}^d\times\mathbb{R}^d$
\begin{eqnarray*}
    &|b(x,y)|\leq L_1\Big(1+|x|+|y|\Big),\\
    &\|\sigma(x,y)\|^2+\int_{\mathbb{Z}}|G(x,y,z)|^2\nu(dz)\leq L_1\Big(1+|x|^2+|y|^2\Big).
\end{eqnarray*}

\noindent{\bf (H3)} the functionals $b$, $\sigma$ and $G$ are non-Lipschitz continuous, and there exist two positive constants $\beta, L_2$ such that for all $(x,\mu)$, $(x',\mu')\in\mathbb{R}^d\times\mathcal{P}_2(\mathbb{R}^d)$ and $z\in\mathbb{Z}$
\begin{eqnarray*}
    &\<x-x',b[x,\mu]-b[x',\mu']\>+\|\sigma[x,\mu]-\sigma[x',\mu']\|^2\leq L_2\Big(\rho(|x-x'|^2)+\rho(\mathcal{W}_2^2(\mu,\mu'))\Big),\\
    &\int_{\mathbb{Z}}|G[x,\mu,z]-G[x',\mu',z]|^2\nu(dz)\leq L_2\Big(\phi(|x-x'|^2)+\phi(\mathcal{W}_2^2(\mu,\mu'))\Big),
\end{eqnarray*}
where $\rho(\cdot)$ and $\phi(\cdot)$ are non-random, non-negative, continuous, concave, non-decreasing functions with $\rho(0)=\phi(0)=0$ and
\begin{eqnarray*}
    \int_{0^+}\Big(\rho(u)+\phi(u)+u\Big)^{-1}du=+\infty.
\end{eqnarray*}

\noindent{\bf (H4)} for every $x\in\overline{\mathcal{D}(A)}$, $x+G(x,y,z)\in\overline{\mathcal{D}(A)}$, $\forall y\in\mathbb{R}^d, z\in\mathbb{Z}$.

\end{hyp}

\begin{remark}\label{rmk2.2}
In fact, Hypothesis \noindent{\bf (H2)} implies that $b$, $\sigma$ and $G$ are continuous in $(x,\mu)$ and satisfy for all $(x,\mu)\in\mathbb{R}^d\times\mathcal{P}(\mathbb{R}^d)$
\begin{eqnarray*}
    |b[x,\mu]|\leq L_1\Big(1+|x|+\|\mu\|_2\Big),
\end{eqnarray*}
\begin{eqnarray*}
    |\sigma[x,\mu]|^2+\int_{\mathbb{Z}}|G[x,\mu,z]|^2\nu(dz)\leq L_1\Big(1+|x|^2+\|\mu\|_2^2\Big),
\end{eqnarray*}
which also implies that
\begin{eqnarray*}
    |b[0,\delta_0]|^2+\|\sigma[0,\delta_0]\|^2<+\infty, \int_{\mathbb{Z}}|G[0,\delta_0,z]|^2\nu(dz)<+\infty.
\end{eqnarray*}
\end{remark}

To address the conditions imposed by Hypothesis \noindent{\bf (H3)}, we need the following fact from Lemma 116 in \cite{Situ}.

\begin{lemma}\label{lem2.4}
Assume that $\varphi(u)$ is a non-negative, increasing function defined on $u\geq0$, $\varphi(0)=0$, $\varphi(u)>0$ for all $u>0$ and $\int_{0^+}1/\varphi(u)du=+\infty$. If for any $t\geq0$, a real non-random function $y_t$ satisfies
\begin{eqnarray*}
    0\leq y_t\leq\int_0^t\varphi(y_s)ds<+\infty,
\end{eqnarray*}
then $y_t=0$ for all $t\geq0$.
\end{lemma}

\subsection{Main Results}

Now we present the main results, whose proofs will be given in the later sections.

The first one is the well-posedness of strong solution under non-Lipschitz conditions and linear growth conditions.

\begin{thm}\label{the3.2}
(Strong solution) If $X_0$ is a $\mathcal{F}_0$-measurable random variable satisfying $\mathbb{E}|X_0|^2<+\infty$ and Hypotheses \noindent{\bf (H1)}-\noindent{\bf (H4)} hold, then there exists a unique strong solution $(X,K)\in\mathcal{A}$ to Eq.(\ref{equ1.1}) and $\mathbb{E}\|X\|_{\mathbb{D}}^2< +\infty$.
\end{thm}

The second one is the existence of weak solution to Eq.(\ref{equ1.1}) under linear growth conditions.

\begin{thm}\label{the4.3}
(Weak solution) Assume that Hypotheses \noindent{\bf (H1)}-\noindent{\bf (H2)} and \noindent{\bf (H4)} hold. Then for any given initial value $X_0$ with finite second moment, Eq.(\ref{equ1.1}) has a weak solution.
\end{thm}

Next, we present the result about the tightness for solutions to a sequence of related equations, which is used to establish the existence of martingale solution. Note that the tightness itself is also of independent interest.

\begin{thm}\label{the5.3}
(Tightness of solutions) For each $n\in \mathbb{N}$, suppose that $(X^n,K^n)$ is a solution to the following equation
\begin{eqnarray*}
    dX^n_t\in-A(X^n_t)dt+b^n[X^n_t,\mathcal{L}_{X^n_t}]dt+\sigma^n[X^n_t,\mathcal{L}_{X^n_t}]dW_t\\
    +\int_{\mathbb{Z}} G^n[X^n_{t-},\mathcal{L}_{X^n_t},z]\widetilde{N}(dt,dz),\ \ t\in[0,1],
\end{eqnarray*}
where the progressively measurable functions $b^n(x,\mu)\rightarrow b(x,\mu)$, $\sigma^n(x,\mu)\rightarrow\sigma(x,\mu)$ and $G^n(x,\mu)\rightarrow G(x,\mu)$ locally in $(x,\mu)$ as $n\rightarrow+\infty$. Assume that $\sup_{n\in \mathbb{N}}\mathbb{E}|X^n_0|^2<+\infty$, then the sequence $(X^n,K^n)_{n\in \mathbb{N}}$ is tight in $\mathbb{D}\times\mathbb{C}$.
\end{thm}

Finally, we establish the existence of martingale solution by the tightness of solutions.

\begin{thm}\label{the5.4}
(Martingale solution) Assume that Hypotheses \noindent{\bf (H1)}-\noindent{\bf (H2)} and \noindent{\bf (H4)} hold. If $\mathcal{L}_{X_0}$ is a probability measure on $\overline{\mathcal{D}(A)}$ such that $\int_{\mathbb{R}^d}|x|^2\mathcal{L}_{X_0}(dx)<+\infty$, then there exists a martingale solution of Eq.(\ref{equ1.1}) with initial distribution $\mathcal{L}_{X_0}$.
\end{thm}

\section{Strong Solution}\label{sec3}
\setcounter{equation}{0}

In this section, we concentrate on proving Theorem \ref{the3.2}, namely the existence and uniqueness of the strong solution to Eq. (\ref{equ1.1}) under Hypotheses \noindent{\bf (H1)}-\noindent{\bf (H4)}. The proof is divided into two steps. First, we establish the existence of the strong solution to Eq. (\ref{equ1.1}) with locally Lipschitz coefficients, building upon the result for globally Lipschitz coefficients. Second, to establish the existence of the strong solution under non-Lipschitz coefficients, we employ a smoothing method to construct a suitable sequence of locally Lipschitz continuous functions that approximate the original coefficients. The uniqueness of the solution is guaranteed by Hypothesis \noindent{\bf (H3)}.

\subsection{Strong Solution under locally Lipschitz Conditions}

Consider the following assumptions:

\noindent{\bf (H3')} The functionals $b$, $\sigma$ and $G$ are globally Lipschitz continuous in the pair $(x,\mu)$, i.e. for all $(x,\mu),(x',\mu')\in\mathbb{R}^d\times\mathcal{P}_2(\mathbb{R}^d)$, there exists some constant $L_3>0$ such that
\begin{eqnarray*}
    &|b[x,\mu]-b[x',\mu']|^2+\|\sigma[x,\mu]-\sigma[x',\mu']\|^2\leq L_3\Big(|x-x'|^2+\mathcal{W}_2^2(\mu,\mu')\Big),\\
    &\int_{\mathbb{Z}}|G[x,\mu,z]-G[x',\mu',z]|^2\nu(dz)\leq L_3\Big(|x-x'|^2+\mathcal{W}_2^2(\mu,\mu')\Big)
\end{eqnarray*}
and
\begin{eqnarray*}
    |b[0,\delta_0]|^2+\|\sigma[0,\delta_0]\|^2<+\infty, \int_{\mathbb{Z}}|G[0,\delta_0,z]|^2\nu(dz)<+\infty.
\end{eqnarray*}

\noindent{\bf (H4')} The functionals $b$, $\sigma$ and $G$ are locally Lipschitz continuous with respect to the state variable and globally Lipschitz continuous on the measure, i.e. for each integer $N>0$, there exist two constants $L_N$ and $L_4$ such that, for all $(x,\mu),(x',\mu')\in\mathbb{R}^d\times\mathcal{P}_2(\mathbb{R}^d)$ with $|x|\vee|x'|\leq N$,
\begin{eqnarray*}
    &|b[x,\mu]-b[x',\mu']|^2+\|\sigma[x,\mu]-\sigma[x',\mu']\|^2\leq L_N|x-x'|^2+L_4\mathcal{W}_2^2(\mu,\mu'),\\
    &\int_{\mathbb{Z}}|G[x,\mu,z]-G[x',\mu',z]|^2\nu(dz)\leq L_N|x-x'|^2+L_4\mathcal{W}_2^2(\mu,\mu').
\end{eqnarray*}

First of all, we study the well-posedness of the solution to Eq.(\ref{equ1.1}) under globally Lipschitz continuous conditions.

\begin{thm}\label{the3.1}
Assume that Hypotheses \noindent{\bf (H1)}, \noindent{\bf (H3')} and \noindent{\bf (H4)} hold. Then for any $\mathcal{F}_0$-measurable random variable $X_0$ such that $\mathbb{E}|X_0|^2<+\infty$, there exists a unique strong solution $(X,K)\in\mathcal{A}$ to Eq.(\ref{equ1.1}).
\end{thm}

\noindent{\bf Proof.}
 We apply the classical Picard iteration method. Let $X^0_t=X_0$ for all $t\in[0,1]$. For each $n=1,2,\cdots$, we can recursively define the pair $(X^n,K^n)$ as the solutions to the following multivalued stochastic differential equation:
\begin{eqnarray*}
    \left\{\begin{array}{ll}
        dX^n_t\in -A(X^n_t)dt+b[X^n_t,\mathcal{L}_{X^{n-1}_t}]dt+\sigma[X^n_t,\mathcal{L}_{X^{n-1}_t}]dW_t\\
        \quad\quad\quad+\int_{\mathbb{Z}}G[X^n_{t-},\mathcal{L}_{X^{n-1}_t},z]\widetilde{N}(dt,dz),~~~~~t\in[0,1],\\
        X^n_0=X_0\in\overline{\mathcal{D}(A)},
    \end{array}\right.
\end{eqnarray*}
where $\mathcal{L}_{X^{n-1}_t}$ denotes the law of $X^{n-1}_t$. More precisely, the pair $(X^n_t,K^n_t)$ satisfies the following recursive equation:
\begin{eqnarray*}
    X^n_t=X_0-K^n_t+\int_0^t b[X^n_s,\mathcal{L}_{X^{n-1}_s}]ds+\int_0^t \sigma[X^n_s,\mathcal{L}_{X^{n-1}_s}]dW_s \nonumber\\
    +\int_0^t \int_{\mathbb{Z}} G[X^n_{s-},\mathcal{L}_{X^{n-1}_s},z]\widetilde{N}(ds,dz),~~\forall t\in[0,1]
\end{eqnarray*}
for $n=1,2,\cdots$.

We first assume that $0\in \text{Int}(\mathcal{D}(A))$ and $0\in A(0)$. Let $\tau_N:=\inf\{t>0:|X_t|\geq N\}$. Applying It\^{o}'s formula, we obtain
\begin{eqnarray*}
    |X^n_{t\wedge\tau_N}|^2&=&|X_0|^2-2\int_0^{t\wedge\tau_N}\langle X^n_s,dK^n_s\rangle+2\int_0^{t\wedge\tau_N}\langle X^n_s,b[X^n_s,\mathcal{L}_{X^{n-1}_s}]\rangle ds\\
    &&+2\int_0^{t\wedge\tau_N}\langle X^n_s,\sigma[X^n_s,\mathcal{L}_{X^{n-1}_s}]dW_s\rangle\\
    &&+2\int_0^{t\wedge\tau_N}\int_{\mathbb{Z}}\langle X^n_s,G[X^n_{s-},\mathcal{L}_{X^{n-1}_s},z]\rangle\widetilde{N}(ds,dz)\\
    &&+\int_0^{t\wedge\tau_N}\|\sigma[X^n_s,\mathcal{L}_{X^{n-1}_s}]\|^2ds+\int_0^{t\wedge\tau_N}\int_{\mathbb{Z}}|G[X^n_{s-},\mathcal{L}_{X^{n-1}_s},z]|^2N(ds,dz)\\
    &=:&|X_0|^2+J_1(t)+\cdots+J_6(t).
\end{eqnarray*}

By Lemma \ref{lem2.1}, we have
\begin{eqnarray*}
    J_1(t)\leq-2\kappa_1\|K^n\|_0^{t\wedge\tau_N}+2\kappa_2\int_0^{t\wedge\tau_N}|X_s^n|ds+2\kappa_3t.
\end{eqnarray*}
Furthermore by Young's inequality, we obtain
\begin{eqnarray*}
    \mathbb{E}\Big(\sup_{r\in[0,t]}J_1(r)\Big)\leq-2\kappa_1\|K^n\|_0^{t\wedge\tau_N}+\kappa_2\int_0^{t\wedge\tau_N}\mathbb{E}\Big(\sup_{r\in[0,s]}|X_r^n|^2\Big)ds+\Big(2\kappa_3+\kappa_2\Big)t.
\end{eqnarray*}

By Hypothesis \noindent{\bf (H3')}, Young's inequality and Remark \ref{rmk2.1}, we get
\begin{eqnarray*}
    J_2(t)&=&2\int_0^{t\wedge\tau_N}\langle X^n_s,b[X^n_s,\mathcal{L}_{X^{n-1}_s}]-b[0,\delta_0]\rangle ds+2\int_0^{t\wedge\tau_N}\langle X^n_s,b[0,\delta_0]\rangle ds\\
    &\leq&2\int_0^{t\wedge\tau_N}|X^n_s|^2ds+\int_0^{t\wedge\tau_N}|b[X^n_s,\mathcal{L}_{X^{n-1}_s}]-b[0,\delta_0]|^2ds+\int_0^{t\wedge\tau_N}|b[0,\delta_0]|^2ds\\
    &\leq&2\int_0^{t\wedge\tau_N}|X^n_s|^2ds+2L_3\int_0^{t\wedge\tau_N}\Big(|X^n_s|^2+\mathcal{W}^2_2(\mathcal{L}_{X^{n-1}_s},\delta_0)\Big)ds\\
    &\leq&(2+2L_3)\int_0^{t\wedge\tau_N}|X^n_s|^2ds+2L_3\int_0^{t\wedge\tau_N}\mathbb{E}(|X^{n-1}_s|^2)ds.
\end{eqnarray*}
Thus
\begin{eqnarray*}
    \mathbb{E}\Big(\sup_{r\in[0,t]}J_2(r)\Big)\leq C\int_0^{t\wedge\tau_N}\mathbb{E}\Big(\sup_{r\in[0,s]}|X^n_r|^2\Big)ds+C\int_0^{t\wedge\tau_N}\mathbb{E}\Big(\sup_{r\in[0,s]}|X^{n-1}_r|^2\Big)ds,
\end{eqnarray*}
where $C$ is a constant that may vary from line to line.

By Burkholder-Davis-Gundy's inequality, Young's inequality, Hypothesis \noindent{\bf (H3')} and Remark \ref{rmk2.1}, we obtain
\begin{eqnarray*}
    \mathbb{E}\Big(\sup_{r\in[0,t]}J_3(r)\Big)&\leq&C\mathbb{E}\langle J_3(\cdot),J_3(\cdot)\rangle _t^{1/2}\\
    &=&C\mathbb{E}\Big(\int_0^{t\wedge\tau_N}|X^n_s|^2\|\sigma[X^n_s,\mathcal{L}_{X^{n-1}_s}]\|^2ds\Big)^{1/2}\\
    &\leq&C\mathbb{E}\Big(\sup_{r\in[0,t]}|X^n_{r\wedge\tau_N}|^2\Big(\int_0^{t\wedge\tau_N}\|\sigma[X^n_s,\mathcal{L}_{X^{n-1}_s}]-\sigma[0,\delta_0]\|^2ds\\
    &&+\int_0^{t\wedge\tau_N}\|\sigma[0,\delta_0]\|^2ds\Big)\Big)^{1/2}\\
    &\leq&C\epsilon\mathbb{E}\Big(\sup_{r\in[0,t]}|X^n_{r\wedge\tau_N}|^2\Big)+C\epsilon^{-1}\int_0^{t\wedge\tau_N}\mathbb{E}\Big(\sup_{r\in[0,s]}|X^n_r|^2\Big)ds\\
    &&+C\epsilon^{-1}\int_0^{t\wedge\tau_N}\|\sigma[0,\delta_0]\|^2ds+C\epsilon^{-1}\int_0^{t\wedge\tau_N}\mathbb{E}\Big(\sup_{r\in[0,s]}|X^{n-1}_r|^2\Big)ds.
\end{eqnarray*}

Similarly, we have
\begin{eqnarray*}
    \mathbb{E}\Big(\sup_{r\in[0,t]}J_4(r)\Big)&\leq&C\mathbb{E}\langle J_4(\cdot),J_4(\cdot)\rangle _t^{1/2}\\
    &\leq&C\epsilon\mathbb{E}\Big(\sup_{r\in[0,t]}|X^n_{r\wedge\tau_N}|^2\Big)+C\epsilon^{-1}\int_0^{t\wedge\tau_N}\mathbb{E}\Big(\sup_{r\in[0,s]}|X^n_r|^2\Big)ds\\
    &&+C\epsilon^{-1}\int_0^{t\wedge\tau_N}\mathbb{E}\Big(\sup_{r\in[0,s]}|X^{n-1}_r|^2\Big)ds\\
    &&+C\epsilon^{-1}\int_0^{t\wedge\tau_N}\int_{\mathbb{Z}}\|G[0,\delta_0,z]\|^2\nu(dz)ds
\end{eqnarray*}
and
\begin{eqnarray*}
    &&\mathbb{E}\Big(\sup_{r\in[0,t]}J_5(r)\Big)+\mathbb{E}\Big(\sup_{r\in[0,t]}J_6(r)\Big)\\
    &\leq& C\int_0^{t\wedge\tau_N}\mathbb{E}\Big(\sup_{r\in[0,s]}|X^n_r|^2\Big)ds+C\int_0^{t\wedge\tau_N}\mathbb{E}\Big(\sup_{r\in[0,s]}|X^{n-1}_r|^2\Big)ds\\
    &&+C\int_0^{t\wedge\tau_N}\|\sigma[0,\delta_0]\|^2ds+C\int_0^{t\wedge\tau_N}\int_{\mathbb{Z}}\|G[0,\delta_0,z]\|^2\nu(dz)ds.
\end{eqnarray*}

Combining the above estimates together, we obtain
\begin{eqnarray*}
    \mathbb{E}\Big(\sup_{r\in[0,t]}|X^n_{r\wedge\tau_N}|^2\Big) &\leq&C+\mathbb{E}|X_0|^2-C\mathbb{E}\|K^n\|_0^{t\wedge\tau_N}+\epsilon\mathbb{E}\Big(\sup_{r\in[0,t]}|X^n_{r\wedge\tau_N}|^2\Big)\\
    &&+C\int_0^{t\wedge\tau_N}\mathbb{E}\Big(\sup_{r\in[0,s]}|X^n_r|^2\Big)ds+C\int_0^{t\wedge\tau_N}\mathbb{E}\Big(\sup_{r\in[0,s]}|X^{n-1}_r|^2\Big)ds.
\end{eqnarray*}
By choosing $\epsilon$ small enough and applying Gronwall's inequality, we get
\begin{eqnarray*}
    \mathbb{E}\Big(\sup_{r\in[0,t]}|X^n_{r\wedge\tau_N}|^2\Big)+C\mathbb{E}\|K^n\|_0^{t\wedge\tau_N}&\leq& C+\mathbb{E}|X_0|^2+C\int_0^{t\wedge\tau_N}\mathbb{E}\Big(\sup_{r\in[0,s]}|X^{n-1}_r|^2\Big)ds.
\end{eqnarray*}
By iteration, we obtain
\begin{eqnarray*}
    \mathbb{E}\Big(\sup_{r\in[0,t]}|X^n_{r\wedge\tau_N}|^2\Big)+C\mathbb{E}\|K^n\|_0^{t\wedge\tau_N}&\leq& C\Big(1+\mathbb{E}|X_0|^2\Big).
\end{eqnarray*}
Let $N\rightarrow+\infty$, by Fatou's lemma, we can deduce that
\begin{eqnarray}\label{equ3.1}
    \sup_{n\in\mathbb{N}}\mathbb{E}\Big(\|X^n\|_{\mathbb{D}}^2\Big)+C\mathbb{E}\|K^n\|_{TV}&\leq& C\Big(1+\mathbb{E}|X_0|^2\Big).
\end{eqnarray}

Let $Z^{n,m}_t:=X^n_t-X^m_t$ and $g_t:=\mathbb{E}\Big(\sup_{s\in[0,t]}|Z^{n,m}_s|^2\Big).$ By It\^{o}'s formula, we have
\begin{eqnarray*}
    |Z^{n,m}_t|^2&:=&-2\int_0^t\langle X^n_s-X^m_s,d(K_s^n-K_s^m)\rangle\\
    &&+2\int_0^t\langle X^n_s-X^m_s,b[X^n_s,\mathcal{L}_{X^{n-1}_s}]-b[X^m_s,\mathcal{L}_{X^{m-1}_s}]\rangle\\
    &&+2\int_0^t\langle X^n_s-X^m_s,(\sigma[X^n_s,\mathcal{L}_{X^{n-1}_s}]-\sigma[X^m_s,\mathcal{L}_{X^{m-1}_s}])dW_s\rangle\\
    &&+2\int_0^t\langle X^n_s-X^m_s,(G[X^n_{s-},\mathcal{L}_{X^{n-1}_s},z]-G[X^m_{s-},\mathcal{L}_{X^{m-1}_s},z])\rangle\widetilde{N}(ds,dz)\\
    &&+\int_0^t\|\sigma[X^n_s,\mathcal{L}_{X^{n-1}_s}]-\sigma[X^m_s,\mathcal{L}_{X^{m-1}_s}]\|^2ds\\ &&+\int_0^t|G[X^n_{s-},\mathcal{L}_{X^{n-1}_s},z]-G[X^m_{s-},\mathcal{L}_{X^{m-1}_s},z]|^2N(ds,dz).
\end{eqnarray*}
Following the similar approach as before, we can derive that
\begin{eqnarray*}
    g_t\leq C\epsilon g(t)+C\epsilon^{-1}\int_0^t\Big(g(s)+\mathcal{W}_2^2(\mathcal{L}_{X^{n-1}_s},\mathcal{L}_{X^{m-1}_s})\Big)ds.
\end{eqnarray*}
By selecting $\epsilon$ sufficiently small such that $C\epsilon<\frac{1}{2}$ and applying Gronwall's lemma, we have
\begin{eqnarray*}
    \mathcal{W}_2^2(\mathcal{L}_{X^n_s},\mathcal{L}_{X^{m}_s})\leq\mathbb{E}\Big(\sup_{s\in[0,t]}|Z^{n,m}_s|^2\Big)\leq C\int_0^t\mathcal{W}_2^2(\mathcal{L}_{X^{n-1}_s},\mathcal{L}_{X^{m-1}_s})ds.
\end{eqnarray*}
Combining this result with (\ref{equ3.1}), through iteration, we derive that there exists a probability measure $\mathcal{L}_{X_t}\in\mathcal{P}(\mathbb{R}^d)$ and a c\`{a}dl\`{a}g $\mathcal{F}_t$-adapted process $X_t$ such that for all $t>0$,
\begin{eqnarray*}
\lim_{n,m\rightarrow+\infty}\mathcal{W}_2^2(\mathcal{L}_{X^n_s},\mathcal{L}_{X^{m}_s})=0
\end{eqnarray*}
and
\begin{eqnarray*}
\lim_{n\rightarrow+\infty}\mathbb{E}\Big(\sup_{t\in[0,1]}|X^n_t-X_t|^2\Big)=0.
\end{eqnarray*}

Define
\begin{eqnarray*}
    K_t&:=&X_0-X_t+\int_0^t b[X_s,\mathcal{L}_{X_s}]ds+\int_0^t \sigma[X_s,\mathcal{L}_{X_s}]dW_s \\
    &&+\int_0^t \int_{\mathbb{Z}} G[X_{s-},\mathcal{L}_{X_s},z]\widetilde{N}(ds,dz),~~\forall t\geq0.
\end{eqnarray*}
It remains to prove that for almost all $\omega$,
\begin{eqnarray}\label{equ3.2}
    \Big(X_{\cdot}(\omega),K_{\cdot}(\omega)\Big)\in\mathcal{A}.
\end{eqnarray}
Note that
\begin{eqnarray*}
    K^n_t&=&X_0-X^n_t+\int_0^t b[X^n_s,\mathcal{L}_{X^n_s}]ds+\int_0^t \sigma[X^n_s,\mathcal{L}_{X^n_s}]dW_s \\
    &&+\int_0^t \int_{\mathbb{Z}} G[X^n_{s-},\mathcal{L}_{X^n_s},z]\widetilde{N}(ds,dz),~~\forall t\geq0.
\end{eqnarray*}
It is easy to prove that
\begin{eqnarray*}
\lim_{n\rightarrow+\infty}\mathbb{E}\Big(\sup_{t\in[0,1]}|K^n_t-K_t|^2\Big)=0.
\end{eqnarray*}
By applying Lemma \ref{lem2.2} and (\ref{equ3.1}), we can obtain (\ref{equ3.2}). The uniqueness of the solution follows by the similar calculations as above.

Next we consider the case that the maximal monotone operator $A$ does not satisfy the assumption that $0\in \text{Int}(\mathcal{D}(A))$ and $0\in A(0)$. This general case can be reduced to a particular one as follows: for any $x\in \text{Int}(\mathcal{D}(A))$, $x^*\in\mathcal{A}(x)$, $z\in\mathbb{Z}$ and any process $\{X_t, t\in[0,1]\}$, define
\begin{eqnarray*}
&&\widetilde{\mathcal{A}}(X_t):=\mathcal{A}(X_t+x)-x^*, \ \ \ \ \ \ \widetilde{b}[X_t,\mathcal{L}_{X_t}]:=b[X_t+x,\mathcal{L}_{X_t-x}]-x^*,\\
&&\widetilde{\sigma}[X_t,\mathcal{L}_{X_t}]:=\sigma[X_t+x,\mathcal{L}_{X_t-x}], \widetilde{G}[X_t,\mathcal{L}_{X_t},z]:=G[X_t+x,\mathcal{L}_{X_t-x},z],\\
&&\widetilde{X}_t:=X_t-x, \widetilde{K}_t:=K_t-tx^*.
\end{eqnarray*}
It is easy to verify that $\widetilde{A}$ is a maximal monotone operator, $0\in \text{Int}(\mathcal{D}(A))$, $0\in\widetilde{A}(0)$ and $\widetilde{b}$, $\widetilde{\sigma}$ and $\widetilde{G}$ satisfy the hypotheses required in the theorem. Furthermore, one can verify that $(\widetilde{X},\widetilde{K})$ is a solution to the following equation:
\begin{eqnarray*}
d\widetilde{X}_t\in -\widetilde{A}(\widetilde{X}_t)dt+\widetilde{b}[\widetilde{X}_t,\mathcal{L}_{\widetilde{X}_t}]dt+\widetilde{\sigma}[\widetilde{X}_t,\mathcal{L}_{\widetilde{X}_t}]dW_t\\
        ~~~~~~~~+\int_{\mathbb{Z}}\widetilde{G}[\widetilde{X}_{t-},\mathcal{L}_{\widetilde{X}_t},z]\widetilde{N}(dt,dz),~~~~t\in[0,1].
\end{eqnarray*}
Finally, it is straightforward to verify through simple calculations that the estimates on $(X,K)$ also hold for $(\tilde{X},\tilde{K})$. $\Box$

Next, we establish the following estimate for the solution of Eq.(\ref{equ1.1}).

\begin{proposition}\label{pro3.1}
Assume that Hypotheses \noindent{\bf (H1)}-\noindent{\bf (H2)} and \noindent{\bf (H4)} hold. If the pair of functions $(X,K)$ is a solution of Eq.(\ref{equ1.1}) and $\mathbb{E}|X_0|^2<+\infty$, then
\begin{eqnarray*}
    \mathbb{E}\Big(\|X\|^2_{\mathbb{D}}+\|K\|_{TV}\Big)\leq C\Big(1+\mathbb{E}|X_0|^2\Big),
\end{eqnarray*}
where $C$ is some positive constant depending on $d$, $L_1$ and $A$.
\end{proposition}
\noindent{\bf Proof.}
Let $\tau_N:=\inf\{t>0:|X_t|>N\}$. For any $t\in[0,1]$, by It\^{o}'s formula we have
\begin{eqnarray*}
    &&|X_{t\wedge\tau_N}|^2\\
    &=&|X_0|^2-2\int^{t\wedge\tau_N}_0\langle X_s,dK_s\rangle+2\int^{t\wedge\tau_N}_0\langle X_s,b[X_s,\mathcal{L}_{X_s}]\rangle ds\\
    &&+2\int^{t\wedge\tau_N}_0\langle X_s,\sigma[X_s,\mathcal{L}_{X_s}]dW_s\rangle +2\int^{t\wedge\tau_N}_0\int_{\mathbb{Z}}\langle X_s,G[X_{s-},\mathcal{L}_{X_s},z]\rangle \widetilde{N}(ds,dz)\\
    &&+\int^{t\wedge\tau_N}_0\|\sigma[X_s,\mathcal{L}_{X_s}]\|^2ds+\int^{t\wedge\tau_N}_0\int_{\mathbb{Z}}|G[X_{s-},\mathcal{L}_{X_s},z]|^2N(ds,dz)\\
    &=:&|X_0|^2+J_1(t)+J_2(t)+\cdots+J_6(t).
\end{eqnarray*}
Next we estimate the above six items separately.

By Hypothesis \noindent{\bf (H2)}, we obtain
\begin{eqnarray*}
    \mathbb{E}\Big(\sup_{r\in[0,t]}J_1(r)\Big)&\leq&-2\kappa_1\mathbb{E}\|K\|_0^{t\wedge\tau_N}+2\kappa_2\mathbb{E}\int_0^{t\wedge\tau_N}|X_s|ds+2\kappa_3t\\
    &\leq&-2\kappa_1\mathbb{E}\|K\|_0^{t\wedge\tau_N}+\kappa_2\int_0^t\mathbb{E}\Big(\sup_{r\in[0,s]}|X_{r\wedge\tau_N}|^2\Big)ds+\Big(2\kappa_3+\kappa_2\Big)t
\end{eqnarray*}
and
\begin{eqnarray*}
    \mathbb{E}\Big(\sup_{r\in[0,t]}\Big(J_2(r)+J_5(r)+J_6(r)\Big)\Big)\leq L_1\int^t_0\Big(1+2\mathbb{E}\Big(\sup_{r\in[0,s]}|X_{r\wedge\tau_N}|^2\Big)\Big)ds.
\end{eqnarray*}
By Hypothesis \noindent{\bf (H2)}, for any $t\in[0,1]$, we have
\begin{eqnarray*}
    \mathbb{E}\int^{t\wedge\tau_N}_0|X_s\cdot\sigma[X_s,\mathcal{L}_{X_s}]|^2ds&=&\mathbb{E}\int^{t\wedge\tau_N-}_0|X_s\cdot\sigma[X_s,\mathcal{L}_{X_s}]|^2ds\\
    &\leq&L_1N^2\int^{t\wedge\tau_N}_0\Big(1+|X_s|^2+\|\mathcal{L}_{X_s}\|^2_2\Big)ds\\
    &\leq&L_1N^2\int^{t\wedge\tau_N}_0\Big(1+N^2+\mathbb{E}|X_s|^2\Big)ds\\
    &\leq&L_1N^2\int^{t\wedge\tau_N}_0\Big(1+2N^2\Big)ds\\
    &<&+\infty.
\end{eqnarray*}
Hence $\Big\{\int^{t\wedge\tau_N}_0\langle X_s,\sigma[X_s,\mathcal{L}_{X_s}]dW_s\rangle\Big\}_{t\in[0,1]}$ is a martingale. Under similar proof as above, we have $\Big\{\int^{t\wedge\tau_N}_0\int_{\mathbb{Z}}\langle X_s,G[X_{s-},\mathcal{L}_{X_s},z]\rangle\widetilde{N}(ds,dz)\Big\}_{t\in[0,1]}$ is also a martingale. Then by Burkholder-Davis-Gundy's inequality and Young's inequality we get
\begin{eqnarray*}
    \mathbb{E}\Big(\sup_{r\in[0,t]}J_3(r)\Big)&\leq& C\mathbb{E}\Big(\int^{t\wedge\tau_N}_0|X_s\cdot\sigma[X_s,\mathcal{L}_{X_s}]|^2ds\Big)^{1/2}\\
    &\leq&C\mathbb{E}\Big(\sup_{r\in[0,t]}|X_{r\wedge\tau_N}|^2\int^{t\wedge\tau_N}_0\|\sigma[X_s,\mathcal{L}_{X_s}]\|^2ds\Big)^{1/2}\\
    &\leq&C\epsilon\mathbb{E}\Big(\sup_{r\in[0,t]}|X_{r\wedge\tau_N}|^2\Big)+C\epsilon^{-1}L_1\int^{t\wedge\tau_N}_0\Big(1+\mathbb{E}|X_s|^2+\|\mathcal{L}_{X_s}\|_2^2\Big)ds\\
    &\leq&C\epsilon\mathbb{E}\Big(\sup_{r\in[0,t]}|X_{r\wedge\tau_N}|^2\Big)+C\epsilon^{-1}L_1\int^{t\wedge\tau_N}_0\Big(1+2\mathbb{E}|X_s|^2\Big)ds\\
    &\leq&C\epsilon\mathbb{E}\Big(\sup_{r\in[0,t]}|X_{r\wedge\tau_N}|^2\Big)+C\epsilon^{-1}L_1\int^t_0\Big(1+2\mathbb{E}\Big(\sup_{r\in[0,s]}|X_{r\wedge\tau_N}|^2\Big)\Big)ds.
\end{eqnarray*}
Similarly,
\begin{eqnarray*}
    \mathbb{E}\Big(\sup_{r\in[0,t]}J_4(r)\Big)\leq C\epsilon\mathbb{E}\Big(\sup_{r\in[0,t]}|X_{r\wedge\tau_N}|^2\Big)+C\epsilon^{-1}L_1\int^t_0\Big(1+2\mathbb{E}\Big(\sup_{r\in[0,s]}|X_{r\wedge\tau_N}|^2\Big)\Big)ds.
\end{eqnarray*}
By choosing $\epsilon$ small enough, we obtain
\begin{eqnarray*}
    \mathbb{E}\Big(\sup_{r\in[0,t]}|X_{r\wedge\tau_N}|^2\Big)\leq \mathbb{E}|X_0|^2 -2\kappa_1\mathbb{E}\|K\|_0^{t\wedge\tau_N}+ C\int^t_0\Big(1+2\mathbb{E}\Big(\sup_{r\in[0,s]}|X_{r\wedge\tau_N}|^2\Big)\Big)ds.
\end{eqnarray*}
By Gronwall's lemma, we have
\begin{eqnarray*}
    \sup_{t\in[0,1]}\Big(\mathbb{E}\Big(\sup_{r\in[0,t]}|X_{r\wedge\tau_N}|^2\Big)+\mathbb{E}\|K\|_0^{t\wedge\tau_N}\Big)\leq C\Big(1+\mathbb{E}|X_0|^2\Big).
\end{eqnarray*}
Let $N\rightarrow+\infty$, by Fatou's lemma, we get the desired result.$\Box$

Finally, we present a generalization of Theorem \ref{the3.1} in which the global Lipschitz condition is replaced by the local Lipschitz condition.

\begin{thm}\label{the3.3}
Assume that Hypotheses \noindent{\bf (H1)}-\noindent{\bf (H2)}, \noindent{\bf (H4')} and \noindent{\bf (H4)} hold. Then for any $\mathcal{F}_0$-adapted random variable $X_0$ such that $\mathbb{E}|X_0|^2<+\infty$, there exists a strong solution $(X,K)\in\mathcal{A}$ to Eq.(\ref{equ1.1}).
\end{thm}

\noindent{\bf Proof.} This theorem is proved by a truncated procedure. For simplicity, we only outline the proof. For each $N>0$, define the truncation function
\begin{eqnarray*}
    b^N[x,\mu]:=\Big\{\begin{array}{ll}
                            b[x,\mu],\ \ \ \ \ \ \ \ \ \ \ if\ |x|\leq N, \\
                            b[Nx/|x|,\mu],\ \ \ if\ |x|>N
                          \end{array}
\end{eqnarray*}
and $\sigma^N$, $G^N$ similarly. Then $b^N$, $\sigma^N$ and $G^N$ satisfy the global Lipschitz condition \noindent{\bf (H3')} and the linear growth condition \noindent{\bf (H2)}. Hence by Theorem \ref{the3.1}, there is a unique strong solution $(X^N,K^N)$ to the equation
\begin{eqnarray}
    X^N_t=X_0-K^N_t+\int_0^t b^N[X^N_s,\mathcal{L}_{X^N_s}]ds+\int_0^t \sigma^N[X^N_s,\mathcal{L}_{X^N_s}]dW_s \nonumber\\
    ~~~~+\int_0^t \int_{\mathbb{Z}} G^N[X^N_{s-},\mathcal{L}_{X^N_s},z]\widetilde{N}(ds,dz),~ t\in[0,1].\label{equ3.8}
\end{eqnarray}
Define the stopping time
\begin{eqnarray*}
    \tau_N=\inf\{t\in[0,1]:|X^N_t|> N\},
\end{eqnarray*}
then it is easy to see that
\begin{eqnarray}\label{equ3.7}
    X^N_t=X^{N+1}_t\leq N,\ \ \ if\ t\in[0,\tau_N].
\end{eqnarray}
This implies $\tau_N\leq\tau_{N+1}$, i.e. $\tau_N$ is increasing. Proposition \ref{pro3.1} shows that $\tau_N\rightarrow1$ as $N\rightarrow+\infty$. By (\ref{equ3.7}), $X_{t\wedge\tau_N}=X^N_{t\wedge\tau_N}$ for $t\in[0,1]$, and it then follows from (\ref{equ3.8}) that
\begin{eqnarray*}
    X_{t\wedge\tau_N}&=&X_0-K_{t\wedge\tau_N}+\int_0^{t\wedge\tau_N} b^N[X_s,\mathcal{L}_{X_s}]ds+\int_0^{t\wedge\tau_N} \sigma^N[X_s,\mathcal{L}_{X_s}]dW_s \\
    &&~~~~+\int_0^{t\wedge\tau_N} \int_{\mathbb{Z}} G^N[X_{s-},\mathcal{L}_{X_s},z]\widetilde{N}(ds,dz)\\
    &=&X_0-K_{t\wedge\tau_N}+\int_0^{t\wedge\tau_N} b[X_s,\mathcal{L}_{X_s}]ds+\int_0^{t\wedge\tau_N} \sigma[X_s,\mathcal{L}_{X_s}]dW_s \\
    &&~~~~+\int_0^{t\wedge\tau_N} \int_{\mathbb{Z}} G[X_{s-},\mathcal{L}_{X_s},z]\widetilde{N}(ds,dz),~ t\in[0,1]\\
\end{eqnarray*}
Letting $N\rightarrow+\infty$ we see that $(X_t,K_t)$ is a solution to Eq.(\ref{equ1.1}).

\subsection{Strong Solution under Non-Lipschitz Conditions}

Next we construct a sequence of Lipschitz continuous  functions to approximate the non-Lipschitz coefficients by smoothing method.

\begin{lemma}\label{lem3.1}
Assume that Hypothesis \noindent{\bf (H2)} holds. Then there exist functions $b^n$, $\sigma^n$ and $G^n$, $n\in\mathbb{N}$ satisfying the following properties:

\begin{enumerate}[(i)]
\item  there exists some positive constant $C>0$ independent of $n$, such that for any $(x,\mu),(x',\mu')\in\mathbb{R}^d\times\mathcal{P}(\mathbb{R}^d)$,
\begin{eqnarray*}
&&|b^n[x,\mu]|\leq C\Big(1+|x|+\|\mu\|_2\Big), \\
&&\|\sigma^n[x,\mu]\|^2+\int_{\mathbb{Z}}|G^n[x,\mu,z]|^2\nu(dz)\leq C\Big(1+|x|^2+\|\mu\|_2^2\Big);
\end{eqnarray*}
\item for each integer $N>0$, there exist some positive constants $C$ and $C_N$ such that for all $|x|\vee|x'|\leq N$,
\begin{eqnarray*}
&&|b^n[x,\mu]-b^n[x',\mu']|\leq k_n\Big(C_N|x-x'|+C\mathcal{W}_2(\mu,\mu')\Big), \\
&&\|\sigma^n[x,\mu]-\sigma^n[x',\mu']\|^2\leq k_n\Big(C_N|x-x'|^2+C\mathcal{W}_2^2(\mu,\mu')\Big),\\
&&\int_{\mathbb{Z}}|G^n[x,\mu,z]-G^n[x',\mu',z]|^2\nu(dz)\leq k_n\Big(C_N|x-x'|^2+C\mathcal{W}_2^2(\mu,\mu')\Big),
\end{eqnarray*}
where $k_n\geq0$ is a constant only dependent on $n$;
\item for any $N>0$,
\begin{eqnarray*}
&&\sup_{|x|\leq N}\sup_{supp(\mu)\subseteq B(0,N)}|b^n[x,\mu]-b[x,\mu]|\rightarrow0,\\
&&\sup_{|x|\leq N}\sup_{supp(\mu)\subseteq B(0,N)}\|\sigma^n[x,\mu]-\sigma[x,\mu]\|^2\rightarrow0,\\
&&\sup_{|x|\leq N}\sup_{supp(\mu)\subseteq B(0,N)}\int_{\mathbb{Z}}|G^n[x,\mu,z]-G[x,\mu,z]|^2\nu(dz)\rightarrow0
\end{eqnarray*}
as $n\rightarrow+\infty$.
\end{enumerate}
\end{lemma}

\begin{remark}
The second conclusion in Lemma \ref{lem3.1} states that the functions $b^n$, $\sigma^n$ and $G^n$ is locally Lipschitz with respect to the state variable and globally Lipschitz on the measure variable.
\end{remark}

\noindent{\bf Proof of Lemma \ref{lem3.1}.}
We apply a smoothing method to the function $b$ with respect to $(x,\mu)$. For any $n\in\mathbb{N}$ and $(x,\mu)\in\mathbb{R}^d\times\mathcal{P}(\mathbb{R}^d)$, define
\begin{eqnarray*}
    &&b^n(x,y):=\iint_{\mathbb{R}^d\times\mathbb{R}^d}b\Big(x-\frac{1}{n}\overline{x},y-\frac{1}{n}\overline{y}\Big)
    J^n(\overline{x})J^n(\overline{y})d\overline{x}d\overline{y},
\end{eqnarray*}
where
\begin{eqnarray*}
    J^n(u)=\left\{\begin{array}{ll}
        a_1\exp\{-\frac{1}{1-|u/b_1|^2}\},~~~~|u|<b_1,\\
        0,~~~~~~~~~~~~~~~~~~~~~~~otherwise
    \end{array}\right.
\end{eqnarray*}
with
\begin{eqnarray*}
    a_0:=\int_{\mathbb{R}^d}\exp\{-\frac{1}{1-|u/b_1|^2}\}du,\ b_1=\Big(\frac{4n}{a_0}\Big)^{1/(d+1)},\ a_1=\frac{b_1}{4n}=(4n)^{-\frac{d}{d+1}}a_0^{-\frac{1}{d+1}}.
\end{eqnarray*}
It is easy to check that $\int_{\mathbb{R}^d}J^n(u)du=1$.

Let $b^n[x,\mu]:=\int_{\mathbb{R}^d}b^n(x,y)\mu(dy)$, then by using Hypothesis \noindent{\bf (H2)}, Fubini's theorem and Young's inequality, we get
\begin{eqnarray*}
    |b^n[x,\mu]|&\leq&\iiint_{\mathbb{R}^d\times\mathbb{R}^d\times\mathbb{R}^d}\Big|b\Big(x-\frac{1}{n}\overline{x},y-\frac{1}{n}\overline{y}\Big)\Big|J^n(\overline{x})J^n(\overline{y})d\overline{x}d\overline{y}\mu(dy)\\
    &\leq&L_1\iiint_{\mathbb{R}^d\times\mathbb{R}^d\times\mathbb{R}^d}\Big(1+\Big|x-\frac{1}{n}\overline{x}\Big|+\Big|y-\frac{1}{n}\overline{y}\Big|\Big)J^n(\overline{x})J^n(\overline{y})d\overline{x}d\overline{y}\mu(dy)\\
    &\leq&L_1\iiint_{\mathbb{R}^d\times\mathbb{R}^d\times\mathbb{R}^d}\Big(1+|x|+|y|+\frac{1}{n}(|\overline{x}|+|\overline{y}|)\Big)J^n(\overline{x})J^n(\overline{y})d\overline{x}d\overline{y}\mu(dy)\\
    &\leq&L_1\Big(\frac{3}{2}+|x|+\frac{1}{2}\|\mu\|_2+\frac{1}{n}\iint_{\mathbb{R}^d\times\mathbb{R}^d}(|\overline{x}|+|\overline{y}|)J^n(\overline{x})J^n(\overline{y})d\overline{x}d\overline{y}\Big)\\
    &\leq&L_1\Big(\frac{3}{2}+|x|+\frac{1}{2}\|\mu\|_2+\frac{2b_1}{n}\Big)\\
    &\leq&C\Big(1+|x|+\|\mu\|_2\Big),
\end{eqnarray*}
which implies the condition (1) for $b^n$.

On the other hand,
\begin{eqnarray*}
    |b^n[x,\mu]-b^n[x',\mu']|\leq|b^n[x,\mu]-b^n[x',\mu]|+|b^n[x',\mu]-b^n[x',\mu']|=:I_1+I_2.
\end{eqnarray*}
By Fubini's lemma, we have
\begin{eqnarray*}
    I_1&=&\Big|\iiint_{\mathbb{R}^d\times\mathbb{R}^d\times\mathbb{R}^d}b\Big(x-\frac{1}{n}\overline{x},y-\frac{1}{n}\overline{y}\Big)J^n(\overline{x})J^n(\overline{y})d\overline{x}d\overline{y}\mu(dy)\\
    &&~-\iiint_{\mathbb{R}^d\times\mathbb{R}^d\times\mathbb{R}^d}b\Big(x'-\frac{1}{n}\overline{x},y-\frac{1}{n}\overline{y}\Big)J^n(\overline{x})J^n(\overline{y})d\overline{x}d\overline{y}\mu(dy)\Big|\\
    &=&n^d\Big|\iiint_{\mathbb{R}^d\times\mathbb{R}^d\times\mathbb{R}^d}b\Big(\tilde{x},y-\frac{1}{n}\overline{y}\Big)J^n(n(x-\tilde{x}))J^n(\overline{y})d\tilde{x}d\overline{y}\mu(dy)\\
    &&~-\iiint_{\mathbb{R}^d\times\mathbb{R}^d\times\mathbb{R}^d}b\Big(\tilde{x},y-\frac{1}{n}\overline{y}\Big)J^n(n(x'-\tilde{x}))J^n(\overline{y})d\tilde{x}d\overline{y}\mu(dy)\Big|\\
    &\leq&n^d\iiint_{\mathbb{R}^d\times\mathbb{R}^d\times\mathbb{R}^d}\Big|b\Big(\tilde{x},y-\frac{1}{n}\overline{y}\Big)\Big|\Big|J^n(n(x-\tilde{x}))-J^n(n(x'-\tilde{x}))\Big|J^n(\overline{y})d\tilde{x}d\overline{y}\mu(dy)\\
    &\leq&n^{d+1}|x-x'|\iiint_{\mathbb{R}^d\times\mathbb{R}^d\times\mathbb{R}^d}\Big(1+|\tilde{x}|+\Big|y-\frac{1}{n}\overline{y}\Big|\Big)\\
    &&~\times\int_0^1\Big|\nabla J^n\Big(n\Big(x-\tilde{x}+\theta(x'-x)\Big)\Big)\Big|d\theta J^n(\overline{y})d\tilde{x}d\overline{y}\mu(dy).
\end{eqnarray*}
It is easy to check that the term $|\nabla J^n(u)|=\frac{2a_1b_1^2|u|}{(b_1^2-u^2)^2}\exp\{-\frac{1}{|u/b_1|^2}\}$ attains its maximal value $\frac{a_1}{b_1}\exp\{-\frac{\sqrt{3}}{\sqrt{3}-1}\}\frac{2/\sqrt[4]{3}}{(1-1/\sqrt{3})^2}\leq\frac{1}{2n}<1$ at $|u|=\frac{b_1}{\sqrt[4]{3}}$. Since $|\overline{y}|\leq b_1$ and $\tilde{x}\in B(x,\frac{b_1}{n})\cap B(x',\frac{b_1}{n})$, we have
\begin{eqnarray*}
    I_1&\leq&n^{d+1}|x-x'|\iiint_{\mathbb{R}^d\times\mathbb{R}^d\times\mathbb{R}^d}\Big(1+|\tilde{x}|+|y|+\frac{1}{n}|\overline{y}|\Big)J^n(\overline{y})d\tilde{x}d\overline{y}\mu(dy)\\
    &\leq&n^{d+1}|x-x'|\Big(C+\|\mu\|_2+\frac{b_1}{n}\Big)\\
    &\leq&Cn^{d+1}|x-x'|,
\end{eqnarray*}
where $C$ may depend on $N$. By Fubini's lemma, we have
\begin{eqnarray*}
    I_2&=&\Big|\iiint_{\mathbb{R}^d\times\mathbb{R}^d\times\mathbb{R}^d}b\Big(x'-\frac{1}{n}\overline{x},y-\frac{1}{n}\overline{y}\Big)J^n(\overline{x})J^n(\overline{y})d\overline{x}d\overline{y}\mu(dy)\\
    &&~-\iiint_{\mathbb{R}^d\times\mathbb{R}^d\times\mathbb{R}^d}b\Big(x'-\frac{1}{n}\overline{x},y-\frac{1}{n}\overline{y}\Big)J^n(\overline{x})J^n(\overline{y})d\overline{x}d\overline{y}\mu'(dy)\Big|\\
    &=&n^d\Big|\iiint_{\mathbb{R}^d\times\mathbb{R}^d\times\mathbb{R}^d}b\Big(x'-\frac{1}{n}\overline{x},\tilde{y}\Big)J^n(\overline{x})J^n(n(y-\tilde{y}))d\overline{x}d\tilde{y}\mu(dy)\\
    &&~-\iiint_{\mathbb{R}^d\times\mathbb{R}^d\times\mathbb{R}^d}b\Big(x'-\frac{1}{n}\overline{x},\tilde{y}\Big)J^n(\overline{x})J^n(n(y-\tilde{y}))d\overline{x}d\tilde{y}\mu'(dy)\Big|\\
    &\leq&n^d\iint_{\mathbb{R}^d\times\mathbb{R}^d}\Big|b\Big(x'-\frac{1}{n}\overline{x},\tilde{y}\Big)\Big|J^n(\overline{x})\Big|\int_{\mathbb{R}^d}J^n(n(y-\tilde{y}))\mu(dy)\\
    &&-\int_{\mathbb{R}^d}J^n(n(y-\tilde{y}))\mu'(dy)\Big|d\overline{x}d\tilde{y}\\
    &\leq&n^d\iint_{\mathbb{R}^d\times\mathbb{R}^d}\Big(1+\Big|x'-\frac{1}{n}\overline{x}\Big|+|\tilde{y}|\Big)\Big|J^n(\overline{x})\Big|\int_{\mathbb{R}^d}J^n(n(y-\tilde{y}))\mu(dy)\\
    &&-\int_{\mathbb{R}^d}J^n(n(y-\tilde{y}))\mu'(dy)\Big|d\overline{x}d\tilde{y}\\
    &\leq&n^d\mathcal{W}_2(\mu,\mu')\iint_{\mathbb{R}^d\times\mathbb{R}^d}\Big(1+|x'|+\frac{1}{n}|\overline{x}|+|\tilde{y}|\Big)\Big|J^n(\overline{x})d\overline{x}d\tilde{y}\\
    &\leq&Cn^d\mathcal{W}_2(\mu,\mu'),
\end{eqnarray*}
where $C$ may depend on $N$ and the third inequality follows by the Kantorovich-
Rubinstein duality relation of 1-order Wasserstein distance
\begin{eqnarray*}
    \mathcal{W}_1(\mu,\mu')=\sup_{f}\Big\{\Big|\int f(y)\mu(dy)-\int f(y)\mu'(dy)\Big|:Lip_f:=\sup_{x\neq y}\frac{|f(x)-f(y)|}{|x-y|}\leq1\Big\},
\end{eqnarray*}
and $\mathcal{W}_1(\mu,\mu')\leq\mathcal{W}_2(\mu,\mu')$ (cf. \cite{MVP}). In fact, as is shown in the proof above, we have $Lip_{J}\leq1$. Thus the condition (2) holds for $b^n$.

By the continuity of $b$, we can obtain that $b$ is uniformly continuous on any compact set, i.e. for any $N>0$ and $\epsilon>0$, there exists a constant $\delta>0$ depending only on $\varepsilon$, such that for $\frac{1}{n}<\delta$ and $|x|\vee|y|\leq N$, $\Big|b(x-\frac{1}{n},y-\frac{1}{n})-b(x,y)\Big|<\epsilon$.
Consequently, for $n\geq\frac{1}{\delta}$,
\begin{eqnarray*}
    &&\sup_{|x|\leq N}\sup_{supp(\mu)\subseteq B(0,N)}\Big|b^n[x,\mu]-b[x,\mu]\Big|\\
    &\leq&\iiint_{\mathbb{R}^d\times\mathbb{R}^d\times\mathbb{R}^d}\sup_{|x|\leq N}\sup_{|y|\leq N}\Big|b(x-\frac{1}{n}\overline{x},y-\frac{1}{n}\overline{y})-b(x,y)\Big|d\overline{x}d\overline{y}\mu(dy)\\
    &<&\epsilon.
\end{eqnarray*}
Thus the condition (3) holds for $b^n$.

Define $\sigma^n$ and $G^n$ similarly as $b^n$. It is easy to prove that $\sigma^n$ and $G^n$ satisfies (1)-(3) by some tedious calculation, so we omit the proofs here.
 $\Box$

Now we are in the position to prove Theorem \ref{the3.2}.

\noindent{\bf Proof of Theorem \ref{the3.2}.}
\noindent{\bf Step 1.} First, we prove the existence of a strong solution. Define $b^n$, $\sigma^n$ and $G^n$ as in Lemma \ref{lem3.1}. In Theorem \ref{the3.3} and Proposition \ref{pro3.1}, we have already established the existence of a family of strong solutions $\{(X^n,K^n)\}_{n\in\mathbb{N}}$ for the following multivalued MVSDEs:
\begin{eqnarray}\label{equ3.3}
    X^n_t=X_0-K^n_t+\int_0^t b^n[X^n_s,\mathcal{L}_{X^n_s}]ds+\int_0^t \sigma^n[X^n_s,\mathcal{L}_{X^n_s}]dW_s \nonumber\\
    ~~~~+\int_0^t \int_{\mathbb{Z}} G^n[X^n_{s-},\mathcal{L}_{X^n_s},z]\widetilde{N}(ds,dz),~~\forall t\in[0,1].
\end{eqnarray}
Moreover,
\begin{eqnarray}\label{equ3.4}
    \sup_{n\in\mathbb{N}}\mathbb{E}\Big(\|X^n\|^2_{\mathbb{D}}+\|K^n\|^2_{TV}\Big)\leq C\Big(1+\mathbb{E}|X_0|^2\Big)\leq k_0<+\infty.
\end{eqnarray}
Consequently, $\lim_{N\rightarrow+\infty}\mathbb{P}(\sup_{n\in\mathbb{N}}\|X^n\|_{\mathbb{D}}\leq N)=1$. Hence for any $t\in[0,1]$,
\begin{eqnarray}\label{equ3.5}
    \mathbb{E}|X^m_t-X^n_t|^2&=&\mathbb{E}|X^m_t-X^n_t|^2I_{\{\|X^m\|_{\mathbb{D}}\leq N,~\|X^n\|_{\mathbb{D}}\leq N\}}\nonumber\\
    &&+\mathbb{E}|X^m_t-X^n_t|^2I_{\{\|X^m\|_{\mathbb{D}}\leq N,~\|X^n\|_{\mathbb{D}}\leq N\}^c}\nonumber\\
    &\leq&\mathbb{E}|X^m_t-X^n_t|^2I_{\{\|X^m\|_{\mathbb{D}}\leq N,~\|X^n\|_{\mathbb{D}}\leq N\}}\nonumber\\
    &&+2k_0\mathbb{E}I_{\{\|X^m\|_{\mathbb{D}}\leq N,~\|X^n\|_{\mathbb{D}}\leq N\}^c},
\end{eqnarray}
where $B^c$ denotes the complement of a set $B$. By  applying It\^{o}'s formula, we have
\begin{eqnarray*}
    &&\mathbb{E}|X^m_t-X^n_t|^2I_{\{\|X^m\|_{\mathbb{D}}\leq N,~\|X^n\|_{\mathbb{D}}\leq N\}}\\
    &\leq&-2\mathbb{E}\int_0^t\langle X^m_s-X^n_s,dK^m_s-dK^n_s\rangle I_{\{\|X^m\|_{\mathbb{D}}\leq N,~\|X^n\|_{\mathbb{D}}\leq N\}}\\
    &&+2\mathbb{E}\int_0^t\langle X^m_s-X^n_s,b^m[X^m_s,\mathcal{L}_{X^m_s}]-b^n[X^n_s,\mathcal{L}_{X^n_s}]\rangle I_{\{\|X^m\|_{\mathbb{D}}\leq N,~\|X^n\|_{\mathbb{D}}\leq N\}}ds\\
    &&+\mathbb{E}\int_0^t\|\sigma^m[X^m_s,\mathcal{L}_{X^m_s}]-\sigma^n[X^n_s,\mathcal{L}_{X^n_s}]\|^2I_{\{\|X^m\|_{\mathbb{D}}\leq N,~\|X^n\|_{\mathbb{D}}\leq N\}}ds\\
    &&+\mathbb{E}\int_0^t\int_{\mathbb{Z}}\|G^m[X^m_{s-},\mathcal{L}_{X^m_s}]-G^n[X^n_{s-},\mathcal{L}_{X^n_s}]\|^2\nu(dz)I_{\{\|X^m\|_{\mathbb{D}}\leq N,~\|X^n\|_{\mathbb{D}}\leq N\}}ds\\
    &=:&J_1(t)+\cdots+J_4(t).
\end{eqnarray*}
By applying Proposition \ref{pro2.1}, we have
\begin{eqnarray*}
    J_1(t)\leq0.
\end{eqnarray*}
For $J_2(t)$, we have
\begin{eqnarray*}
    J_{2}(t)&\leq&2\mathbb{E}\int_0^t\<X^m_s-X^n_s,b^m[X^m_s,\mathcal{L}_{X^m_s}]-b[X^m_s,\mathcal{L}_{X^m_s}]\>I_{\{\|X^m\|_{\mathbb{D}}\leq N,~\|X^n\|_{\mathbb{D}}\leq N\}}ds\\
    &&+2\mathbb{E}\int_0^t\<X^m_s-X^n_s,b[X^m_s,\mathcal{L}_{X^m_s}]-b[X^n_s,\mathcal{L}_{X^n_s}]\>I_{\{\|X^m\|_{\mathbb{D}}\leq N,~\|X^n\|_{\mathbb{D}}\leq N\}}ds\\
    &&+2\mathbb{E}\int_0^t\<X^m_s-X^n_s,b[X^n_s,\mathcal{L}_{X^n_s}]-b^n[X^n_s,\mathcal{L}_{X^n_s}]\>I_{\{\|X^m\|_{\mathbb{D}}\leq N,~\|X^n\|_{\mathbb{D}}\leq N\}}ds\\
    &=:&I_1(t)+I_2(t)+I_3(t).
\end{eqnarray*}
Applying condition (3) in Lemma \ref{lem3.1}, we have
\begin{eqnarray*}
    \lim_{m\rightarrow+\infty}\Big|b^m[X^m_s,\mathcal{L}_{X^m_s}]-b[X^m_s,\mathcal{L}_{X^m_s}]\Big|I_{\{\|X^m\|_{\mathbb{D}}\leq N\}}=0.
\end{eqnarray*}
By applying dominated convergence theorem, we obtain
\begin{eqnarray*}
    \limsup_{m,n\rightarrow+\infty}(I_1(t)+I_3(t))\leq2\limsup_{m,n\rightarrow+\infty}\mathbb{E}\int_0^t|X^m_s-X^n_s|^2I_{\{\|X^m\|_{\mathbb{D}}\leq N,~\|X^n\|_{\mathbb{D}}\leq N\}}ds.
\end{eqnarray*}
By Hypothesis \noindent{\bf (H3)}, we have
\begin{eqnarray*}
    I_2(t)\leq2L_2\mathbb{E}\int_0^t\Big(\rho\left(|X^m_s-X^n_s|^2\right)+\rho\left(\mathcal{W}_2^2(\mathcal{L}_{X^m_s},\mathcal{L}_{X^n_s})\right)\Big)ds.
\end{eqnarray*}
Hence
\begin{eqnarray*}
    \limsup_{m,n\rightarrow+\infty}J_2(t)&\leq&2\limsup_{m,n\rightarrow+\infty}\mathbb{E}\int_0^t|X^m_s-X^n_s|^2I_{\{\|X^m\|_{\mathbb{D}}\leq N,~\|X^n\|_{\mathbb{D}}\leq N\}}ds\\
    &&+2L_2\limsup_{m,n\rightarrow+\infty}\mathbb{E}\int_0^t\Big(\rho\left(|X^m_s-X^n_s|^2\right)+\rho\left(\mathcal{W}_2^2(\mathcal{L}_{X^m_s},\mathcal{L}_{X^n_s})\right)\Big)ds\\
    &\leq&2\int_0^t\limsup_{m,n\rightarrow+\infty}\mathbb{E}|X^m_s-X^n_s|^2I_{\{\|X^m\|_{\mathbb{D}}\leq N,~\|X^n\|_{\mathbb{D}}\leq N\}}ds\\
    &&4L_2\int_0^t\rho\Big(\limsup_{m,n\rightarrow+\infty}\mathbb{E}|X^m_s-X^n_s|^2\Big)ds.
\end{eqnarray*}
Similar results also hold for $J_3(t)$ and $J_4(t)$. Therefore
\begin{eqnarray*}
    &&\limsup_{m,n\rightarrow+\infty}\mathbb{E}|X^m_t-X^n_t|^2\\
    &\leq&2k_0\mathbb{E}I_{\{\|X^m\|_{\mathbb{D}}\leq N,~\|X^n\|_{\mathbb{D}}\leq N\}^c}\\
    &&+2\int_0^t\limsup_{m,n\rightarrow+\infty}\mathbb{E}|X^m_s-X^n_s|^2I_{\{\|X^m\|_{\mathbb{D}}\leq N,~\|X^n\|_{\mathbb{D}}\leq N\}}ds\\
    &&+6L_2\int_0^t\rho\Big(\limsup_{m,n\rightarrow+\infty}\mathbb{E}|X^m_s-X^n_s|^2\Big)ds\\
    &&+2L_2\int_0^t\phi\Big(\limsup_{m,n\rightarrow+\infty}\mathbb{E}|X^m_s-X^n_s|^2\Big)ds.
\end{eqnarray*}
Letting $N\rightarrow+\infty$, we can obtain that
\begin{eqnarray}\label{equ3.9}
    \limsup_{m,n\rightarrow+\infty}\mathbb{E}|X^m_t-X^n_t|^2\leq C\int_0^t\psi\Big(\limsup_{m,n\rightarrow+\infty}\mathbb{E}|X^m_s-X^n_s|^2\Big)ds,
\end{eqnarray}
where $\psi(u)=u+\rho(u)+\phi(u)$ is a non-random, non-negative, continuous, concave and non-decreasing function such that $\psi(0)=0$ and $\int_{0^+}du/\psi(u)=+\infty$. By applying Lemma \ref{lem2.4}, we have
\begin{eqnarray*}
    \limsup_{m,n\rightarrow+\infty}\mathbb{E}|X^m_t-X^n_t|^2=0.
\end{eqnarray*}
Hence there exists a $\mathcal{F}_t$-adapted process $X_t$, such that
\begin{eqnarray*}
    \limsup_{n\rightarrow+\infty}\mathbb{E}|X^n_t-X_t|^2=0.
\end{eqnarray*}
So $X^n_t\rightarrow X_t$, in probability for each $t$, and one can choose a subsequence $\{n_k\}$ of $\{n\}$, denoted by $\{n\}$ again, such that $\mathbb{P}$-a.s. as $n\rightarrow+\infty$,
\begin{eqnarray*}
X^n_t\rightarrow X_t, \ \forall\  t\in[0,1].
\end{eqnarray*}
By Fatou's lemma, we have
\begin{eqnarray*}
    \mathbb{E}(\sup_{t\in[0,1]}|X_t|^2)&=&\mathbb{E}(\sup_{t\in[0,1]}\lim_{n\rightarrow+\infty}|X^n_t|^2)\leq\mathbb{E}(\liminf_{n\rightarrow+\infty}\sup_{t\in[0,1]}|X^n_t|^2)\\
    &\leq&\liminf_{n\rightarrow+\infty}\mathbb{E}(\sup_{t\in[0,1]}|X^n_t|^2)\leq k_0.
\end{eqnarray*}
For any $t\in[0,1]$, define
\begin{eqnarray*}
    K_t&:=&X_0-X_t+\int_0^t b[X_s,\mathcal{L}_{X_s}])ds+\int_0^t \sigma[X_s,\mathcal{L}_{X_s}]dW_s \\
    &&+\int_0^t \int_{\mathbb{Z}} G[X_{s-},\mathcal{L}_{X_s},z]\widetilde{N}(ds,dz),
\end{eqnarray*}
where $\mathcal{L}_{X_s}=\mathbb{P}\circ X_s^{-1}$. It remains to prove that for almost all $\omega$,
\begin{eqnarray}\label{equ3.6}
    \Big(X_{\cdot}(\omega),K_{\cdot}(\omega)\Big)\in\mathcal{A}.
\end{eqnarray}
Note that
\begin{eqnarray*}
    K^n_t&=&X_0-X^n_t+\int_0^t b^n[X^n_s,\mathcal{L}_{X^n_s}]ds+\int_0^t \sigma^n[X^n_s,\mathcal{L}_{X^n_s}]dW_s \\
    &&+\int_0^t \int_{\mathbb{Z}} G^n[X^n_{s-},\mathcal{L}_{X^n_s},z]\widetilde{N}(ds,dz).
\end{eqnarray*}
It is easy to prove that for any $t\in[0,1]$,
\begin{eqnarray*}
\lim_{n\rightarrow+\infty}\mathbb{E}|K^n_t-K_t|^2=0.
\end{eqnarray*}
By using Lemma \ref{lem2.2} and (\ref{equ3.4}), we get (\ref{equ3.6}) and therefore the existence.

\noindent{\bf Step 2.} The pathwise uniqueness of the strong solution is guaranteed by Hypothesis \noindent{\bf (H3)}. The proof, which follows from similar calculations as (\ref{equ3.9}), is omitted here for brevity. $\Box$

\section{Weak Solution}\label{sec4}
\setcounter{equation}{0}

In this section we discuss the existence of weak solution to Eq.(\ref{equ1.1}). We prove that for a sequence of solutions $(X^n_t,K^n_t,W_t,\zeta_t)$, $n=1,2,\cdots$, to Eq.(\ref{equ4.2}), by Skorohod's weak convergence technique, one can construct a convergent sequence of processes $(\widetilde{X}^{n_k}_t,\widetilde{K}^{n_k}_t,\widetilde{W}^{n_k}_t,\widetilde{\zeta}^{n_k}_t)$, $k=1,2,\cdots$, where $\{n_k\}$ is a subsequence of $\{n\}$, such that for each $k=1,2,\cdots$, $(\widetilde{X}^{n_k}_t,\widetilde{K}^{n_k}_t,\widetilde{W}^{n_k}_t,\widetilde{\zeta}^{n_k}_t)$ possesses the same finite-dimensional probability distributions as those of $(X^{n_k}_t,K^{n_k}_t,W_t,\zeta_t)$. Subsequently, the existence of a weak solution can be established through the limit process of the sequence $(\widetilde{X}^{n_k}_t,\widetilde{K}^{n_k}_t,\widetilde{W}^{n_k}_t,\widetilde{\zeta}^{n_k}_t)$, $k=1,2,\cdots$.

\subsection{Preliminary Lemma}
Define $b^n$, $\sigma^n$ and $G^n$, $n=1,2,\cdots$, as in the proof of Lemma \ref{lem3.1}. Then for each $n=1,2,\cdots$, there exists a unique strong solution $(X^n_t,K^n_t)$ to the following SDE:
\begin{eqnarray}
    X^n_t&=&X_0-K^n_t+\int_0^t b^n[X^n_s,\mathcal{L}_{X^n_s}]ds+\int_0^t \sigma^n[X^n_s,\mathcal{L}_{X^n_s}]dW_s \nonumber\\
    &&+\int_0^t \int_{\mathbb{Z}} G^n[X^n_{s-},\mathcal{L}_{X^n_s},z]\widetilde{N}(ds,dz),\ \ t\in[0,1].\label{equ4.2}
\end{eqnarray}
Denote $q(dt,dz)=\widetilde{N}(dt,dz)$ the Poisson martingale measure with compensator $\nu(dz)dt$ such that $q(dt,dz)=p(dt,dz)-\nu(dz)dt$ and $p(dt,dz)=N(dt,dz)$. Define
\begin{eqnarray}\label{equ4.9}
    \zeta_t:=\int_0^t\int_{|z|\leq1}zq(ds,dz)+\int_0^t\int_{|z|>1}zp(ds,dz).
\end{eqnarray}
Obviously, $\zeta_t$ is a compound Poisson process with the compensated Poisson measure $q(dt,dz)$.

\begin{lemma}\label{lem4.1}
Assume that for every $n=1,2,\cdots$, $(X^n_t,K^n_t)$ is the solution of Eq.(\ref{equ4.2}) driven by $(W_t,\zeta_t)$, then there exists a sequence of c\`{a}dl\`{a}g processes $(\widetilde{X}^{n_k}_t,\widetilde{K}^{n_k}_t,\widetilde{W}^{n_k}_t,\widetilde{\zeta}^{n_k}_t)$, with $n_0=0$ and $\{n_k\}$ a subsequence of $\{n\}$, defined on some probability space ($\widetilde{\Omega},\widetilde{\mathcal{F}},\widetilde{\mathbb{P}}$) such that
\begin{enumerate}[(i)]
\item all finite-dimensional probability distributions of $(\widetilde{X}^{n_k}_t,\widetilde{K}^{n_k}_t,\widetilde{W}^{n_k}_t,\widetilde{\zeta}^{n_k}_t)$ coincide with the finite-dimensional probability distributions of $(X^{n_k}_t,K^{n_k}_t,W_t,\zeta_t)$, $k=1,2,\cdots$;
\item $\widetilde{\eta}^{n_k}_t\rightarrow\widetilde{\eta}^0_t$, in probability, as $k\rightarrow+\infty$, $\forall t\in[0,1]$, as $\widetilde{\eta}^{n_k}_t=\widetilde{X}^{n_k}_t,\widetilde{K}^{n_k}_t,\widetilde{W}^{n_k}_t,\widetilde{\zeta}^{n_k}_t$, $k=1,2,\cdots$;
\item for all $k=1,2,\cdots$, $(\widetilde{X}^{n_k}_t,\widetilde{K}^{n_k}_t)$ satisfies the following SDE:
\begin{eqnarray*}
    \widetilde{X}^{n_k}_t=X_0+\widetilde{K}^{n_k}_t+\int_0^t b^{n_k}[\widetilde{X}^{n_k}_s,\mathcal{L}_{\widetilde{X}^{n_k}_s}]ds+\int_0^t \sigma^{n_k}[\widetilde{X}^{n_k}_s,\mathcal{L}_{\widetilde{X}^{n_k}_s}]d\widetilde{W}^{n_k}_s,\\
    +\int_0^t \int_{\mathbb{Z}} G^{n_k}[\widetilde{X}^{n_k}_{s-},\mathcal{L}_{\widetilde{X}^{n_k}_s},z]\widetilde{q}^{n_k}(ds,dz),\ \ t\in[0,1],
\end{eqnarray*}
where $\widetilde{q}^{n_k}(ds,dz)$ is a compensated Poisson measure corresponding to $\widetilde{\zeta}^{n_k}_t$ in the form of (\ref{equ4.9}) for all $k=0,1,2,\cdots$. Moreover, for all $k=0,1,2,\cdots$, $\widetilde{q}^{n_k}(dt,dz)$ are Poisson martingale measures with the same compensator $\nu(dz)dt$, $\widetilde{W}^{n_k}_t$ are Brownian motions on the probability space $(\widetilde{\Omega},\widetilde{\mathcal{F}},\widetilde{\mathbb{P}})$.
\end{enumerate}
\end{lemma}

\noindent{\bf Proof.}
Applying Lemma \ref{lem3.1} and Proposition \ref{pro3.1}, we can immediately deduce that
\begin{eqnarray}\label{equ4.1}
    \sup_{n\in\mathbb{N}}\mathbb{E}\Big(\|X^n\|^2_{\mathbb{D}}+\|K^n\|^2_{TV}\Big)\leq C\Big(1+\mathbb{E}|X_0|^2\Big)\leq k_0<+\infty.
\end{eqnarray}
Hence for the sequence $(X^n_t,K^n_t,W_t,\zeta_t)$, $n=1,2,\cdots$, by Skorohod's theorem \cite{Situ}, we only need to prove that
\begin{eqnarray}\label{equ4.3}
\Big\{\begin{array}{ll}
     \lim_{N\rightarrow+\infty}\sup_{n\in\mathbb{N}}\sup_{0\leq t\leq1}\mathbb{P}(|\eta^n_t|>N)=0, \\
     \lim_{h\downarrow0}\sup_{n\in\mathbb{N}}\sup_{|t-r|\leq h}\mathbb{P}(|\eta^n_t-\eta^n_r|>\epsilon)=0.
\end{array}
\end{eqnarray}
as $\eta^n_t=X^n_t,K^n_t,W^n_t,\zeta^n_t$, $n=1,2,\cdots$.

For any $0<r\leq t\leq1$, $n\geq1$, we have
\begin{eqnarray*}
    \mathbb{E}|X^n_t-X^n_r|^2&=&-2\mathbb{E}\int_r^t\langle X^n_s-X^n_r,dK^n_s-dK^n_r\rangle\nonumber\\
    &&+2\mathbb{E}\int_r^t\langle X^n_s-X^n_r,b^n[X^n_s,\mathcal{L}_{X^n_s}]\rangle ds\nonumber\\
    &&+2\mathbb{E}\int_r^t\|\sigma^n[X^n_s,\mathcal{L}_{X^n_s}]\|^2ds\nonumber\\
    &&+2\mathbb{E}\int_r^t\int_{\mathbb{Z}}|G^n[X^n_{s-},\mathcal{L}_{X^n_s},z]|^2\nu(dz)ds\nonumber\\
    &\leq&\mathbb{E}\int_r^t|X^n_s-X^n_r|^2ds+C\mathbb{E}\int_r^t\Big(1+|X^n_s|^2+\|\mathcal{L}_{X^n_s}\|_2^2\Big)ds\nonumber\\
    &\leq&4k_0(t-r)+C(1+2k_0)(t-r)\nonumber\\
    &\leq&C(t-r).
\end{eqnarray*}
Thus for $h>0$,
\begin{eqnarray*}
    \sup_{n\in\mathbb{N}}\sup_{|t-r|\leq h}\mathbb{E}|X^n_t-X^n_r|^2\leq Ch,
\end{eqnarray*}
which implies that for any $0\leq t,r\leq1$ and $\epsilon>0$,
\begin{eqnarray*}
    \lim_{h\downarrow0}\sup_{n\in\mathbb{N}}\sup_{|t-r|\leq h}\mathbb{P}(|X^n_t-X^n_r|>\epsilon)\leq\lim_{h\downarrow0}\frac{Ch}{\epsilon^2}=0.
\end{eqnarray*}
Moreover, by (\ref{equ4.1}) we have
\begin{eqnarray*}
    \lim_{N\rightarrow+\infty}\sup_{n\in\mathbb{N}}\sup_{0\leq t\leq1}\mathbb{P}(|X^n_t|>N)\leq\lim_{N\rightarrow+\infty}\frac{k_0}{N^2}=0.
\end{eqnarray*}
Therefore (\ref{equ4.3}) holds for $X^n_t$. Denote
\begin{eqnarray*}
    \zeta_t=\int_0^t\int_{|z|\leq1}z\widetilde{N}(ds,dz)+\int_0^t\int_{|z|>1}zN(ds,dz)=:\zeta^1_t+\zeta^2_t.
\end{eqnarray*}
From the proof of Lemma 173 in \cite{Situ}, we know that $\zeta_t$ and $W_t$ also satisfy (\ref{equ4.3}).

Finally, we shall prove that $K^n_t$ satisfies (\ref{equ4.3}). In fact,
\begin{eqnarray*}
    K^n_t-K^n_r&=&X^n_r-X^n_t\int_r^t b^n[X^n_s,\mathcal{L}_{X^n_s}]ds+\int_r^t \sigma^n[X^n_s,\mathcal{L}_{X^n_s}]dW_s \\
    &&+\int_r^t \int_{\mathbb{Z}} G^n[X^n_{s-},\mathcal{L}_{X^n_s},z]q(ds,dz).
\end{eqnarray*}
Following the above method, we obtain
\begin{eqnarray*}
    \mathbb{E}|K^n_t-K^n_r|^2\leq C(t-r).
\end{eqnarray*}
Combining the above estimate with (\ref{equ4.1}), we obtain (\ref{equ4.3}) for $K^n_t$, which completes the proof. $\Box$

\begin{remark}\label{rmk4.1}
By this Lemma, if we can prove that as $k\rightarrow+\infty$,
\begin{eqnarray}
    &\Big|\int_0^t\Big(b^{n_k}[\widetilde{X}^{n_k}_s,\mathcal{L}_{\widetilde{X}^{n_k}_s}]-b[\widetilde{X}^0_s,\mathcal{L}_{\widetilde{X}^0_s}]\Big)ds\Big|\rightarrow0,~in~\widetilde{\mathbb{P}},\label{equ4.6}\\
    &\int_0^t\sigma^{n_k}[\widetilde{X}^{n_k}_s,\mathcal{L}_{\widetilde{X}^{n_k}_s}]d\widetilde{W}^{n_k}_s\rightarrow\int_0^t \sigma[\widetilde{X}^0_s,\mathcal{L}_{\widetilde{X}^0_s}]d\widetilde{W}^0_s,~in~\widetilde{\mathbb{P}},\label{equ4.7}\\
    &\int_0^t \int_{\mathbb{Z}} G^{n_k}[\widetilde{X}^{n_k}_{s-},\mathcal{L}_{\widetilde{X}^{n_k}_s},z]\widetilde{q}^{n_k}(ds,dz)\rightarrow\int_0^t \int_{\mathbb{Z}} G[\widetilde{X}^0_{s-},\mathcal{L}_{\widetilde{X}^0_s},z]\widetilde{q}^0(ds,dz),~in~\widetilde{\mathbb{P}},\label{equ4.8}
\end{eqnarray}
then $(\widetilde{X}^0_t,\widetilde{K}^0_t)$ is a weak solution to Eq.(\ref{equ1.1}).
\end{remark}

\subsection{Existence of Weak Solutions}

We proceed to prove the existence of a weak solution to Eq.(\ref{equ1.1}).

\noindent{\bf Proof of Theorem \ref{the4.3}.}
By Lemma \ref{lem3.1} we can approximate $b$, $\sigma$ and $G$ smoothly by $b^n$, $\sigma^n$ and $G^n$, respectively. Then there exists a sequence of pathwise unique strong solution $(X^n,K^n)$ to Eq.(\ref{equ4.2}) driven by $(W,\zeta)$, where $\zeta$ is defined as (\ref{equ4.9}). To prove the existence of a weak solution, with Lemma \ref{lem4.1} we need to show that (\ref{equ4.6})-(\ref{equ4.8}) in Remark \ref{rmk4.1} hold. Next we only give the proofs of (\ref{equ4.8}) since the proofs for (\ref{equ4.6}) and (\ref{equ4.7}) are similar.

By the proof of Theorem \ref{the3.2}, we have
\begin{eqnarray*}
    \sup_{n\in\mathbb{N}}\mathbb{E}\|\widetilde{X}^n\|^2_{\mathbb{D}}+\sup_{n\in\mathbb{N}}\mathbb{E}\|\widetilde{K}^n\|^2_{TV}\leq k_0<+\infty.
\end{eqnarray*}
By Lemma \ref{lem4.1}, there exists a subsequence $\{n_k\}$ of $\{n\}$, $k=1,2,\cdots$, and $n_0=0$, such that $\widetilde{X}^{n_k}_t\rightarrow \widetilde{X}^0_t$, $\widetilde{K}^{n_k}_t\rightarrow \widetilde{K}^0_t$, in probability for each $t\in[0,1]$, as $k\rightarrow+\infty$. We can then choose a subsequence $\{n_k^l\}$ of $\{n_k\}$, $l=1,2,\cdots$, denoted by $\{n\}$ again, such that $\widetilde{X}^n_t\rightarrow \widetilde{X}^0_t$, $\widetilde{K}^n_t\rightarrow \widetilde{K}^0_t$, $\mathbb{P}-a.s.$ for all $t\in[0,1]$ as $n\rightarrow+\infty$.
Hence by Fatou's lemma, we have
\begin{eqnarray*}
    \mathbb{E}\|\widetilde{X}^0\|^2_{\mathbb{D}}=\mathbb{E}\Big(\sup_{t\in[0,1]}\lim_{n\rightarrow+\infty}|\widetilde{X}^n_t|^2\Big)\leq\liminf_{n\rightarrow+\infty}\mathbb{E}\|\widetilde{X}^n\|^2_{\mathbb{D}}\leq k_0<+\infty.
\end{eqnarray*}
Similarly,
\begin{eqnarray*}
    \mathbb{E}\|\widetilde{K}^0\|^2_{TV}\leq k_0<+\infty.
\end{eqnarray*}

For brevity, we denote $\{n_k\}$ by $\{n\}$ again. By Remark 397 in \cite{Situ}, we may assume that for all $t\in[0,1]$, $n=1,2,\cdots$, $\widetilde{X}^n_{t}$ and $\widetilde{X}^0_{t}$ are uniformly bounded, that is, there exists some constant $\bar{k}_0>0$ such that $\sup_{t\in[0,1]}|\widetilde{X}_t^n|\leq \overline{k}_0$ and $\sup_{t\in[0,1]}|\widetilde{X}_t^0|\leq \overline{k}_0$ for any $n=1,2,\cdots$. 

Now for any given $\epsilon>0$,
\begin{eqnarray*}
    &&\widetilde{\mathbb{P}}\Big(\Big|\int_0^t\int_{\mathbb{Z}}G^n[\widetilde{X}^n_{s-},\mathcal{L}_{\widetilde{X}^n_s},z]\widetilde{q}^n(ds,dz)-\int_0^t\int_{\mathbb{Z}}G[\widetilde{X}^0_{s-},\mathcal{L}_{\widetilde{X}^0_s},z]\widetilde{q}^0(ds,dz)\Big|>\epsilon\Big)\\
    &\leq&\widetilde{\mathbb{P}}\Big(\Big|\int_0^t\int_{\mathbb{Z}}\Big(G^n[\widetilde{X}^n_{s-},\mathcal{L}_{\widetilde{X}^n_s},z]-G[\widetilde{X}^n_{s-},\mathcal{L}_{\widetilde{X}^n_s},z]\Big)\widetilde{q}^n(ds,dz)\Big|>\epsilon/3\Big)\\
    &&+\widetilde{\mathbb{P}}\Big(\Big|\int_0^t\int_{\mathbb{Z}}\Big(G[\widetilde{X}^n_{s-},\mathcal{L}_{\widetilde{X}^n_s},z]-G[\widetilde{X}^0_{s-},\mathcal{L}_{\widetilde{X}^0_s},z]\Big)\widetilde{q}^n(ds,dz)\Big|>\epsilon/3\Big)\\
    &&+\widetilde{\mathbb{P}}\Big(\Big|\int_0^t\int_{\mathbb{Z}}G[\widetilde{X}^0_{s-},\mathcal{L}_{\widetilde{X}^0_s},z]\widetilde{q}^n(ds,dz)-\int_0^t\int_{\mathbb{Z}}G[\widetilde{X}^0_{s-},\mathcal{L}_{\widetilde{X}^0_s},z]\widetilde{q}^0(ds,dz)\Big|>\epsilon/3\Big)\\
    &=:& I^n_1+I^n_2+I^n_3.
\end{eqnarray*}

By Chebyshev's inequality we have
\begin{eqnarray*}
    I^n_1\leq\frac{9}{\epsilon^2}\mathbb{E}^{\widetilde{\mathbb{P}}}\int_0^t\int_{\mathbb{Z}}\Big|G^n[X,\mathcal{L}_X,z]-G[X,\mathcal{L},z]\Big|^2\nu(dz)ds.
\end{eqnarray*}
Then by condition (3) in Lemma \ref{lem3.1}, we get that $\lim_{n\rightarrow+\infty}I^n_1=0$.

Also by Chebyshev's inequality we obtain as $n\rightarrow+\infty$,
\begin{eqnarray}\label{equ4.10}
    I^n_2\leq\frac{9}{\epsilon^2}\mathbb{E}^{\widetilde{\mathbb{P}}}\int_0^t\int_{\mathbb{Z}}\Big|G[\widetilde{X}^n_{s-},\mathcal{L}_{\widetilde{X}^n_s},z]-G[\widetilde{X}^0_{s-},\mathcal{L}_{\widetilde{X}^0_s},z]\Big|^2\nu(dz)ds\rightarrow0.
\end{eqnarray}
In fact for each $s\in[0,t]$, $\widetilde{X}^n_s\rightarrow\widetilde{X}^0_s$ in probability implies that $\mathcal{L}_{\widetilde{X}^n_s}\rightarrow\mathcal{L}_{\widetilde{X}^0_s}$ weakly. Combining the above estimate with the assumptions: $\sup_{s\in[0,1]}|\widetilde{X}_s^n|\leq \overline{k}_0$ for any $n=1,2,\cdots$ and $\sup_{s\in[0,1]}|\widetilde{X}_s^0|\leq \overline{k}_0$, we can deduce that
\begin{eqnarray*}
    \lim_{n\rightarrow+\infty}\int_{\mathbb{R}^d}|x|^2\mathcal{L}_{\widetilde{X}^n_s}(dx)&=&\lim_{n\rightarrow+\infty}\int_{\mathbb{R}^d}|x|^2I_{\{|x|\leq \overline{k}_0\}}\mathcal{L}_{\widetilde{X}^n_s}(dx)\\
    &=&\int_{\mathbb{R}^d}|x|^2I_{\{|x|\leq \overline{k}_0\}}\mathcal{L}_{\widetilde{X}^0_s}(dx)\\
    &=&\int_{\mathbb{R}^d}|x|^2\mathcal{L}_{\widetilde{X}^0_s}(dx)
\end{eqnarray*}
where the second equality holds since the integrand is actually a bounded continuous function. Hence by the Definition 6.8 in Villani \cite{Villani1} we have
$\mathcal{W}_2(\mathcal{L}_{\widetilde{X}^n_s},\mathcal{L}_{\widetilde{X}^0_s})\rightarrow0$ for any $s\in[0,1]$ as $n\rightarrow+\infty$. Since continuous function $G$ can preserve the convergence in probability, we have
\begin{eqnarray*}
    \int_{\mathbb{Z}}\Big|G[\widetilde{X}^n_{s-},\mathcal{L}_{\widetilde{X}^n_s},z]-G[\widetilde{X}^0_{s-},\mathcal{L}_{\widetilde{X}^0_s},z]\Big|^2\nu(dz)\rightarrow0,\ in\ \widetilde{\mathbb{P}}.
\end{eqnarray*}
Applying Lebesgue's dominated convergence theorem, we conclude that (\ref{equ4.10}) holds.

As stated in Remark 397 in \cite{Situ}, it is noted that
\begin{eqnarray*}
    I^n_3&\leq&2\Big(\frac{12}{\epsilon}\Big)^2\mathbb{E}\int_0^1\int_{0<|z|<\delta}|G[\widetilde{X}^0_{s-},\mathcal{L}_{\widetilde{X}^0_s},z]|^2\nu(dz)ds\\
    &&+\widetilde{\mathbb{P}}\Big(\Big|\int_0^t\int_{|z|\geq\delta}G[\widetilde{X}^0_{s-},\mathcal{L}_{\widetilde{X}^0_s},z]\widetilde{q}^n(ds,dz)\\
    &&\ \ \ \ \ \ \ \ -\int_0^t\int_{|z|\geq\delta}G[\widetilde{X}^0_{s-},\mathcal{L}_{\widetilde{X}^0_s},z]\widetilde{q}^0(ds,dz)\Big|>\epsilon/6\Big)\\
    &=:&J_1^{0,\delta}+J_2^{0,\delta}.
\end{eqnarray*}

We can choose $\delta>0$ sufficiently small such that $J_1^{0,\delta}<\widetilde{\epsilon}/3$. On the other hand, for any $\epsilon>0$, we also have
\begin{eqnarray*}
    J_2^{0,\delta}&\leq&2\Big(\frac{18}{\epsilon}\Big)^2\mathbb{E}\int_0^1\int_{|z|\geq\delta}I_{\{\|\widetilde{X}^0\|_{\mathbb{D}}\leq k_0\}}\Big|G[\widetilde{X}^0_{s-},\mathcal{L}_{\widetilde{X}^0_s},z]\\
    &&~~~~~~~-\sum^{2^m-1}_{i=0}G[\widetilde{X}^0_{i/2^m-},\mathcal{L}_{\widetilde{X}^0_{i/2^m}},z]I_{(\frac{i}{2^m},\frac{i+1}{2^m}]}(s)\Big|^2\nu(dz)ds\\
    &&+\widetilde{\mathbb{P}}\Big(\sum_{i=0}^{2^m-1}I_{\{\|\widetilde{X}^0\|_{\mathbb{D}}\leq k_0\}}\Big|\int_{i/2^m}^{(i+1)/2^m}\int_{|z|\geq\delta}G[\widetilde{X}^0_{i/2^m-},\mathcal{L}_{\widetilde{X}^0_{i/2^m}},z]\widetilde{q}^n(ds,dz)\\
    &&~~~~~~~-\int_{i/2^m}^{(i+1)/2^m}\int_{|z|\geq\delta}G[\widetilde{X}^0_{i/2^m-},\mathcal{L}_{\widetilde{X}^0_{i/2^m}},z]\widetilde{q}^0(ds,dz)\Big|>\epsilon/6\Big)\\
    &=&J^m_1+J^m_2,
\end{eqnarray*}
where $0<\frac{1}{2^m}<\frac{2}{2^m}<\cdots<\frac{i}{2^m}<\cdots<1$ is a partition of $[0,1]$. By Hypothesis \noindent{\bf (H2)} we have $\lim_{m\rightarrow+\infty}J^m_1=0$. Hence we can choose $m$ large enough such that $J^m_1<\frac{\widetilde{\epsilon}}{6}$. Note that $G[\widetilde{X}^0_{i/2^m-},\mathcal{L}_{\widetilde{X}^0_{i/2^m}},z]$ is a finite Borel measurable function of $z$ for any $i$ and $m$, then by Lemma 400 in \cite{Situ}, for fixed $m$, $\delta$, there exists some constant $\bar{N}$ such that for any $n\geq\bar{N}$, $J^m_2<\frac{\widetilde{\epsilon}}{6}$. Consequently, we deduce that $\lim_{n\rightarrow+\infty}I^n_3=0$. This implies that $\lim_{n\rightarrow+\infty}\sum_{i=1}^3I^n_i=0$, which verifies Eq.(\ref{equ4.8}) in Remark \ref{rmk4.1}. Thus $(\widetilde{X}^0_t,\widetilde{K}^0_t)$ is a weak solution of Eq(\ref{equ1.1}). $\Box$

\section{Martingale Solution}\label{sec5}
\setcounter{equation}{0}

This section is devoted to establish the existence of martingale solution to Eq.(\ref{equ1.1}) by the Aldous' tightness criterion.

\subsection{Tightness}
First we will show that any class of solutions of Eq.(\ref{equ1.1}) is $tight$. To this purpose, we introduce the following equivalent characterization of tightness (cf. \cite{PB}).

\begin{thm}\label{the5.1} (Aldous' tightness criterion)
Let $\{X^n\}_{n\in \mathbb{N}}$ be a sequence of real-valued stochastic processes taking values in $\mathbb{D}$. The family $\{\mathcal{L}_{X^n}\}_{n\in \mathbb{N}}$ is tight if the following conditions hold:
\begin{enumerate}[(i)]
\item uniform boundedness:
\begin{eqnarray}\label{equ5.2}
    \lim_{a\rightarrow+\infty}\sup_{n\in \mathbb{N}}\mathbb{P}\{\|X^n\|_{\mathbb{D}}\geq a\}=0;
\end{eqnarray}
\item Aldous' condition: for any $\rho>0$,
\begin{eqnarray}\label{equ5.3}
    \lim_{\delta\downarrow0}\sup_{n\in \mathbb{N}}\sup_{(S,S')\in A_{\delta}}\mathbb{P}\{|X^n_{S'}-X^n_S|> \rho\}=0,
\end{eqnarray}
where $A_\delta$ is the set of all pairs of stopping times $(S,S')$ such that $0\leq S\leq S'\leq S+\delta\leq1$ a.s..
\end{enumerate}
\end{thm}

\noindent{\bf Proof of Theorem \ref{the5.3}.}
As in the proof of Theorem \ref{the3.1}, we assume that $0\in \text{Int}(\mathcal{D}(A))$ and $0\in A(0)$.

Define $b^n$, $\sigma^n$ and $G^n$ as in Lemma \ref{lem3.1}. By Theorem \ref{the3.3}, for each $n\in\mathbb{N}$, there exists a unique strong solution $(X^n,K^n)$ for the following multivalued MVSDE:
\begin{eqnarray*}
    X^n_t=X^n_0-K^n_t+\int_0^t b^n[X^n_s,\mathcal{L}_{X^n_s}]ds+\int_0^t \sigma^n[X^n_s,\mathcal{L}_{X^n_s}]dW_s \nonumber\\
    ~~~~+\int_0^t \int_{\mathbb{Z}} G^n[X^n_{s-},\mathcal{L}_{X^n_s},z]\widetilde{N}(ds,dz),~~\forall t\in[0,1].
\end{eqnarray*}
By Aldous' tightness criteria, we only need to verify that (\ref{equ5.2})-(\ref{equ5.3}) hold for $(X^n,K^n)$.

By Proposition \ref{pro3.1},
\begin{eqnarray}\label{equ5.4}
    \sup_{n\in\mathbb{N}}\mathbb{E}\Big(\|X^n\|^2_{\mathbb{D}}+\|K^n\|^2_{TV}\Big)\leq C\Big(1+\sup_{n\in\mathbb{N}}\mathbb{E}|X^n_0|^2\Big)\leq k_0<+\infty.
\end{eqnarray}
Consequently,
\begin{eqnarray*}
    \lim_{a\rightarrow+\infty}\sup_{n\in \mathbb{N}}\mathbb{P}(\|X^n\|_{\mathbb{D}}\geq a)\leq\lim_{a\uparrow+\infty}\sup_{n\in \mathbb{N}}\frac{\mathbb{E}\|X^n\|^2_{\mathbb{D}}}{a^2}\leq\lim_{a\uparrow+\infty}\frac{k_0}{a^2}=0,
\end{eqnarray*}
which implies  (\ref{equ5.2}).

Next, it remains to verify  (\ref{equ5.3}). For any given $\rho>0$ and $a>0$, we have
\begin{eqnarray}
    &&\mathbb{P}\Big(|X^n_{S'}-X^n_S|>\rho\Big)\nonumber\\
    &\leq&\mathbb{P}\Big(|X^n_{S'}-X^n_S|I_{\{\|X^n\|_{\mathbb{D}}\leq a\}}>\rho\Big)+\mathbb{P}\Big(\|X^n\|_{\mathbb{D}}>a\Big)\nonumber\\
    &\leq&\frac{1}{\rho^2}\mathbb{E}\Big(|X^n_{S'}-X^n_S|^2I_{\{\|X^n\|_{\mathbb{D}}\leq a\}}\Big)+\frac{1}{a^2}\mathbb{E}\|X^n\|^2_{\mathbb{D}}.\label{equ5.5}
\end{eqnarray}

For any $\epsilon>0$ and $a>0$, let
\begin{eqnarray*}
    \mathbb{K}_{\epsilon,a}:=\{x\in\mathbb{R}^d:|x-y|\geq\epsilon, \forall~y\notin\overline{\mathcal{D}(A)}\}\cap\overline{B(0,a)}.
\end{eqnarray*}
 Then $\mathbb{K}_{\epsilon,a}$ is a convex compact set contained in $\text{Int}(\mathcal{D}(A))$, and it is nonempty if $\epsilon$ is smaller than a certain $\epsilon_0>0$ (independent of $a$).

Next let $f_a: (0,+\infty)\rightarrow[0,+\infty)$ be the function defined by
\begin{eqnarray*}
    f_a(\epsilon):=\sup\{|x^*|: x^*\in A(x),\ x\in\mathbb{K}_{\epsilon,a}\}.
\end{eqnarray*}
(By convention, $\sup\emptyset=0$). Since $A$ is locally bounded on $\text{Int}(\mathcal{D}(A))$, this function is well-defined and non-increasing. For $\delta>0$ and $a>0$, define
\begin{eqnarray*}
    g_a(\delta):=\inf\{\epsilon\in(0,\epsilon_0): f_a(\epsilon)\leq\delta^{-1/2}\}.
\end{eqnarray*}
(By convention, $\inf\emptyset=\epsilon_0$). We have $f_a(\delta+g_a(\delta))\leq\delta^{-1/2}$ for every $\delta>0$ and $\lim_{\delta\downarrow0}g_a(\delta)=0$. Let $\delta_a>0$ such that $\delta_a+g_a(\delta_a)<\epsilon_0$.

For a pair of fixed but arbitrary stopping time $(S,S')$,
\begin{eqnarray*}
    X^n_{S'}-X^n_S&=&K^n_{S}-K^n_{S'}+\int_{S}^{S'}b^n[X^n_r,\mathcal{L}_{X^n_r}]dr+\int_{S}^{S'}\sigma^n[X^n_r,\mathcal{L}_{X^n_r}]dW^n_r\\
    &&+\int_{S}^{S'}\int_{\mathbb{Z}}G^n[X^n_{r-},\mathcal{L}_{X^n_r},z]\widetilde{N}(dr,dz).
\end{eqnarray*}
Thus by It\^{o}'s formula,
\begin{eqnarray*}
    |X^n_{S'}-X^n_S|^2&=&2\int_{S}^{S'}\langle X^n_S-X^n_r,dK^n_r\rangle+2\int_{S}^{S'}\langle X^n_r-X^n_S,b^n[X^n_r,\mathcal{L}_{X^n_r}]\rangle dr\\
    &&+2\int_{S}^{S'}\langle X^n_r-X^n_S,\sigma^n[X^n_r,\mathcal{L}_{X^n_r}]dW^n_r\rangle\\
    &&+2\int_{S}^{S'}\int_{\mathbb{Z}}\langle X^n_r-X^n_S,G^n[X^n_{r-},\mathcal{L}_{X^n_r},z]\rangle \widetilde{N}(dr,dz)\\
    &&+\int_{S}^{S'}\|\sigma^n[X^n_r,\mathcal{L}_{X^n_r}]\|^2dr\\
    &&+\int_{S}^{S'}\int_{\mathbb{Z}}|G^n[X^n_{r-},\mathcal{L}_{X^n_r},z]|^2N(dr,dz)\\
    &=:&J_1+J_2+\cdots+J_6.
\end{eqnarray*}
Fix $a>0$ and $\delta\in(0,\delta_a\wedge1]$. Then $\mathbb{K}_{\delta+g_a(\delta),a}$ is nonempty. Hence there exists a projection of $X^n_S$ on $\mathbb{K}_{\delta+g_a(\delta),a}$, denoted as $Y^{n,\delta,a}_S$, moreover
\begin{eqnarray*}
    |X^n_S-Y^{n,\delta,a}_S|\leq\delta+g_a(\delta),~~\mathbb{P}\text{-}a.s.~on~\{\|X^n\|_{\mathbb{D}}\leq a\}
\end{eqnarray*}
since $X^n_S\in\overline{\mathcal{D}(A)}$. Let $z=Y^{n,\delta,a}_S$, $z^*\in A(z)$ ($A(Y^{n,\delta,a}_S)$ is nonempty), we obtain that
\begin{eqnarray*}
    J_1&\leq&2\int_{S}^{S'}\langle X^n_S-Y^{n,\delta,a}_S,dK^n_r\rangle+2\int_{S}^{S'}\langle Y^{n,\delta,a}_S-X^n_r,dK^n_r-z^*dr\rangle\\
    &&+2\int_{S}^{S'}\langle Y^{n,\delta,a}_{\tau_n}-X^n_r,z^*\rangle dr\\
    &\leq&2(\delta+g_a(\delta))\|K^n\|_{TV}+4a\delta f_a(\delta+g_a(\delta))\\
    &\leq&2(\delta+g_a(\delta))\|K^n\|_{TV}+4\delta^{1/2}a,
\end{eqnarray*}
$\mathbb{P}$-$a.s.$ on $\{\|X^n\|_{\mathbb{D}}\leq a\}$.

By Remark \ref{rmk2.2} and Young's inequality, we can deduce that
\begin{eqnarray*}
    J_2&\leq&4a\int_{S}^{S'}|b^n[X^n_r,\mathcal{L}_{X^n_r}]|dr\\
    &\leq&4aC\int_{S}^{S'}\Big(1+|X^n_r|+\|\mathcal{L}_{X^n_r}\|_2\Big)dr\\
    &\leq&4aC\int_{S}^{S'}\Big(1+a+(\mathbb{E}|X^n_r|^2)^{1/2}\Big)dr\\
    &\leq&4aC\delta(1+a).
\end{eqnarray*}
Similarly,
\begin{eqnarray*}
    J_5+J_6&\leq&C\int_{S}^{S'}\Big(1+|X^n_r|^2+\|\mathcal{L}_{X^n_r}\|^2_2\Big)dr\\
    &\leq&C\int_{S}^{S'}\Big(1+a^2+\mathbb{E}|X^n_r|^2\Big)dr\\
    &\leq&C\delta(1+a^2).
\end{eqnarray*}
Hence
\begin{eqnarray}
    \mathbb{E}\Big(|X^n_{S'}-X^n_S|^2I_{\{\|X^n\|_{\mathbb{D}}\leq a\}}\Big) &\leq&2(\delta+g_a(\delta))\|K^n\|_{TV}+4\delta^{1/2}a\nonumber\\
    &&+4aC\delta(1+a)+C\delta(1+a^2).\label{equ5.6}
\end{eqnarray}
Combining $\lim_{\delta\downarrow0}g_a(\delta)=0$ with (\ref{equ5.5})-(\ref{equ5.6}), we can obtain that (\ref{equ5.3}) holds for $X^n$. Hence $\{X^n\}_{n\in \mathbb{N}}$ is tight in $\mathbb{D}$.

For any $t\in[0,1]$,
\begin{eqnarray*}
    K^n_t=X_0-X^n_t+\int_0^t b^n[X^n_s,\mathcal{L}_{X^n_s}]ds+\int_0^t \sigma^n[X^n_s,\mathcal{L}_{X^n_s}]dW_s \\
    ~~~~+\int_0^t \int_{\mathbb{Z}} G^n[X^n_{s-},\mathcal{L}_{X^n_s},z]\widetilde{N}(ds,dz),~~\mathbb{P}^n\text{-}a.s..
\end{eqnarray*}
Note that
$\{\int_0^t b^n[X^n_s,\mathcal{L}_{X^n_s}]ds\}$, $\{\int_0^t \sigma^n[X^n_s,\mathcal{L}_{X^n_s}]dW_s\}$ and $\{\int_0^t \int_{\mathbb{Z}} G^n[X^n_{s-},\mathcal{L}_{X^n_s},z]\widetilde{N}(ds,dz)\}$ are all tight in $\mathbb{D}$, hence $\{K^n\}_{n\in \mathbb{N}}$ is also tight in $\mathbb{D}$, which completes the proof. $\Box$

\subsection{Existence of Martingale Solution}
We begin by introducing an essential set $F_a:=\{(x,\eta)\in\bar{\Omega}_0:\|\eta\|_{TV}\leq a\}$, $a\geq0$. Obviously $F_a$ is closed in $\bar{\Omega}$.

\noindent{\bf Proof of Theorem \ref{the5.4}.}
Define $b^n$, $\sigma^n$ and $G^n$ as in the proof of Lemma \ref{lem3.1}. We can then find a stochastic basis $(\Omega,\mathcal{F},\mathbb{P},\{\mathcal{F}_t\}_{t\geq0})$, a standard, $d$-dimensional Brownian motion $W$, a compensated Poisson measure $\widetilde{N}$ with respect to this basis and a $\mathcal{F}_0$-measurable $d$-dimensional random vector $X_0$ such that $\mathcal{L}_{X_0}=\mathbb{P}\circ X_0^{-1}$.

By Theorem \ref{the3.3} and Proposition \ref{pro3.1}, there exists a family of pathwise unique strong solution $\{(X^n,K^n)\}_{n\in\mathbb{N}}$ for the following multivalued MVSDEs:
\begin{eqnarray*}
    \Bigg\{\begin{array}{lll}
        dX^n_t\in A(X^n_t)dt+b^n[X^n_t,\mathcal{L}_{X^n_t}]dt+\sigma^n[X^n_t,\mathcal{L}_{X^n_t}]dW_t\\
        \quad\quad\quad+\int_{\mathbb{Z}}G^n[X^n_{t-},\mathcal{L}_{X^n_t},z]\widetilde{N}(dt,dz),~~t\in[0,1],\\
        X^n_0=X_0.
    \end{array}
\end{eqnarray*}

Let $R_n:=\mathbb{P}\circ(X^n,K^n)^{-1}$. By Theorem \ref{the5.3} and Prohorov's theorem, the sequence $\{R_n\}_{n\in\mathbb{N}}$ is relatively compact in $\mathcal{P}(\bar{\Omega})$. So there exists a subsequence (also denoted by $\{R_n\}$) converging to some probability measure $R\in\mathcal{P}(\bar{\Omega})$. It is left to prove that $R$ is the solution of the martingale problem for Eq.(\ref{equ1.1}).

Observe that $R_n(\bar{\Omega}_0)=1$. The canonical processes $\bar{X}^n:\bar{\Omega}\rightarrow\mathbb{D}$ and $\bar{K}^n:\bar{\Omega}\rightarrow\mathbb{C}$, which are extended from $(X^n,K^n)$ via (\ref{equ5.8}), are progressively measurable with respect to the stochastic basis $(\bar{\Omega},\bar{\mathcal{F}}^{R_n},R_n,\{\bar{\mathcal{F}}^{R_n}_t\}_{t\geq0})$. By Theorem \ref{the3.3}, there exists a constant $C>0$ such that
\begin{eqnarray*}
    \sup_{n\in\mathbb{N}}\mathbb{E}^{R_n}\Big(\|\bar{X}^n\|^2_{\mathbb{D}}+\sup_{n\in\mathbb{N}}\|\bar{K}^n\|^2_{TV}\Big)\leq C\Big(1+\mathbb{E}|X_0|^2\Big)\leq k_0<+\infty.
\end{eqnarray*}
We begin by prove that $R(\bar{\Omega}_0)=1$. Indeed, since the set $F_a$ is closed for each $a\geq0$, by Chebyshev's inequality we have
\begin{eqnarray*}
    R(\bar{\Omega}_0)&\geq& R(F_a)\geq\limsup_{n\rightarrow+\infty}R_n(F_a)=1-\liminf_{n\rightarrow+\infty}R_n(\|\bar{K}^n\|_{TV}>a)\\
    &\geq&1-\frac{1}{a^2}\liminf_{n\rightarrow+\infty}\mathbb{E}^{R_n}\|\bar{K}^n\|^2_{TV}\geq1-\frac{C}{a^2}\Big(1+\|\mathcal{L}_{X_0}\|_2^2\Big).
\end{eqnarray*}
Let $n\in\mathbb{N}$ and $f\in C^2_c(\mathbb{R}^d)$, define
\begin{eqnarray*}
    \bar{M}_t^{f,n}&:=&f(\bar{X}^n_t)-f(\bar{X}_0)+\int_0^t \langle \nabla f(\bar{X}^n_s), d\bar{K}^n_s\rangle-\int_0^t\mathfrak{L}^n(\mathcal{L}_{\bar{X}^n_s})f(\bar{X}^n_s)ds\\
    &&-\frac{1}{2}\int_0^t\int_{\mathbb{Z}}{G^n}^*[\bar{X}^n_{s-},\mathcal{L}_{\bar{X}^n_s},z]\nabla^2f(\bar{X}^n_s)G^n[\bar{X}^n_{s-},\mathcal{L}_{\bar{X}^n_s},z]\nu(dz)ds,~t\in[0,1],
\end{eqnarray*}
where $\mathcal{L}_{\bar{X}^n_s}=R_n\circ (\bar{X}^n_s)^{-1}$ and $\mathfrak{L}^n(\mu)f(x):=\frac{1}{2}tr\Big(\sigma^n\sigma^{n*}[x,\mu]\nabla^2\Big)f(x)+\langle b^n[x,\mu],\nabla f(x)\rangle$ for $(x,\mu)\in\mathbb{D}\times\mathcal{P}(\mathbb{D})$.

For a probability measure $R$ on $\bar{\Omega}$, an integrable $\bar{\mathcal{F}}^R_t$-adapted process $M$ is a martingale under $(\bar{\Omega},\bar{\mathcal{F}}^R,R,\{\bar{\mathcal{F}}^R_t\}_{t\geq0})$ if and only if $\mathbb{E}^R(M_t-M_s)\Gamma=0$ holds for each $0\leq s<t\leq1$ and for every bounded, real-valued, continuous, $\bar{\mathcal{F}}_s$-measurable function $\Gamma$.

Fix $f\in C^2_c(\mathbb{R}^d)$, by It\^{o}'s formula, we can deduce that
\begin{eqnarray*}
    \bar{M}_t^{f,n}=\int_0^t\langle\nabla f(\bar{X}^n_s),\sigma^n[\bar{X}^n_s,\mathcal{L}_{\bar{X}^n_s}]dW_s\rangle+\int_0^t\int_{\mathbb{Z}}\langle\nabla f(\bar{X}^n_s),G^n[\bar{X}^n_{s-},\mathcal{L}_{\bar{X}^n_s},z]\rangle\widetilde{N}(ds,dz)
\end{eqnarray*}
is a martingale under $R_n$ due to $\sup_{n\in\mathbb{N}}\mathbb{E}^{R_n}\|\bar{X}^n\|^2_{\mathbb{D}}\leq k_0<+\infty$. Consequently, for any $0\leq s<t\leq1$,
\begin{eqnarray}\label{equ5.9}
    \mathbb{E}^{R_n}(\bar{M}_t^{f,n}-\bar{M}_s^{f,n})\Gamma=0, ~~\forall~n\in\mathbb{N}.
\end{eqnarray}
It remains to prove that we can pass the limit as $n \rightarrow+\infty$ in (\ref{equ5.9}).

Since $R_n\rightarrow R$, by Skorohod's representation theorem, we can construct  a probability space $(\bar{\Omega},\bar{\mathcal{F}},\bar{\mathbb{P}})$, a family of random variables $\{(\bar{X}^n,\bar{K}^n)\}_{n\in\mathbb{N}}$ and $(\bar{X},\bar{K})$ defined on this probability space, such that
\begin{enumerate}[(i)]
\item $(\bar{X}^n,\bar{K}^n)$ (resp. $(\bar{X},\bar{K})$ has the law $R_n$ (resp. $R$), for each $n\in\mathbb{N}$;
\item $(\bar{X}^n,\bar{K}^n)\rightarrow(\bar{X},\bar{K})$ in $\bar{\Omega}$, $\bar{\mathbb{P}}$-a.s., as $n\rightarrow+\infty$.
\end{enumerate}

Naturally, (ii) implies that $(\bar{X}^n,\bar{K}^n)\rightarrow(\bar{X},\bar{K})$ in distribution. Note that the function $(f(\bar{X}^n_t)-f(\bar{X}^n_s))\Gamma$ is bounded and continuous within the uniform convergence topology. Hence by Portmanteau's theorem,
\begin{eqnarray}\label{equ5.10}
    \lim_{n\rightarrow+\infty}\mathbb{E}^{R_n}\Big(f(\bar{X}^n_t)-f(\bar{X}^n_s)\Big)\Gamma=\mathbb{E}^{R}\Big(f(\bar{X}_t)-f(\bar{X}_s)\Big)\Gamma.
\end{eqnarray}
On the other hand, since $\langle \nabla f(\bar{X}^n_s), b^n[\bar{X}^n_s,\mathcal{L}_{\bar{X}^n_s}]\rangle\rightarrow\langle \nabla f(\bar{X}_s), b[\bar{X}_s,\mathcal{L}_{\bar{X}_s}]\rangle$ in distribution uniformly on compact subsets of $\bar{\Omega}$ for almost all $s\in[0,1]$, we obtain
\begin{eqnarray*}
    \lim_{n\rightarrow+\infty}\mathbb{E}^{R_n}\langle \nabla f(\bar{X}^n_s), b^n[\bar{X}^n_s,\mathcal{L}_{\bar{X}^n_s}]\rangle\Gamma=\mathbb{E}^{R}\langle \nabla f(\bar{X}_s), b[\bar{X}_s,\mathcal{L}_{\bar{X}_s}]\rangle\Gamma.
\end{eqnarray*}
Similarly, we get
\begin{eqnarray*}
    &&\lim_{n\rightarrow+\infty}\mathbb{E}^{R_n}\Big(\sum_{i,j=1}^d \Big(\nabla^2f(\bar{X}^n_s)\Big)_{ij}\Big({\sigma^n}^*[\bar{X}^n_s,\mathcal{L}_{\bar{X}^n_s}]\sigma^n[\bar{X}^n_s,\mathcal{L}_{\bar{X}^n_s}]\Big)_{ij}\Big)\Gamma\\
    &=&\mathbb{E}^{R}\Big(\sum_{i,j=1}^d \Big(\nabla^2f(\bar{X}_s)\Big)_{ij}\Big({\sigma}^*[\bar{X}_s,\mathcal{L}_{\bar{X}_s}]\sigma[\bar{X}_s,\mathcal{L}_{\bar{X}_s}]\Big)_{ij}\Big)\Gamma
\end{eqnarray*}
and
\begin{eqnarray*}
    &&\lim_{n\rightarrow+\infty}\mathbb{E}^{R_n}\Big(\int_{\mathbb{Z}}{G^n}^*[\bar{X}^n_{s-},\mathcal{L}_{\bar{X}^n_s},z]\nabla^2f(\bar{X}^n_s)G^n[\bar{X}^n_{s-},\mathcal{L}_{\bar{X}^n_s},z]\nu(dz)\Big)\Gamma\\
    &=&\mathbb{E}^{R}\Big(\int_{\mathbb{Z}}G^*[\bar{X}_{s-},\mathcal{L}_{\bar{X}_s},z]\nabla^2f(\bar{X}_s)G[\bar{X}_{s-},\mathcal{L}_{\bar{X}_s},z]\nu(dz)\Big)\Gamma.
\end{eqnarray*}
By (\ref{equ5.10}), Lebesgue's dominated convergence theorem, and Fubini's theorem, we have for almost all $0\leq s<t\leq1$,
\begin{eqnarray}\label{equ5.11}
    &&\lim_{n\rightarrow+\infty}\mathbb{E}^{R_n}\Big(f(\bar{X}^n_t)-f(\bar{X}^n_s)-\int_s^t \mathfrak{L}^n(\mathcal{L}_{\bar{X}^n_s})f(\bar{X}^n_s)ds\nonumber\\
    &&\ \ \ -\frac{1}{2}\int_s^t\int_{\mathbb{Z}}{G^n}^*[\bar{X}^n_{s-},\mathcal{L}_{\bar{X}^n_s},z]\nabla^2f(\bar{X}^n_s)G^n[\bar{X}^n_{s-},\mathcal{L}_{\bar{X}^n_s},z]\nu(dz)ds\Big)\Gamma\nonumber\\
    &=&\mathbb{E}^{R}\Big(f(\bar{X}_t)-f(\bar{X}_s)-\int_s^t \mathfrak{L}(\mathcal{L}_{\bar{X}^n_s})f(\bar{X}^n_s)ds\nonumber\\
    &&\ \ \ -\frac{1}{2}\int_s^t\int_{\mathbb{Z}}G^*[\bar{X}_{s-},\mathcal{L}_{\bar{X}_s},z]\nabla^2f(\bar{X}_s)G[\bar{X}_{s-},\mathcal{L}_{\bar{X}_s},z]\nu(dz)ds\Big)\Gamma.
\end{eqnarray}
By Proposition 3.4 in \cite{Zalinescu}, we can deduce that
\begin{eqnarray}\label{equ5.12}
    \lim_{n\rightarrow+\infty}\mathbb{E}^{R_n}\Big(\Gamma\int_0^t \langle \nabla f(\bar{X}^n_s), d\bar{K}^n_s\rangle\Big)=\mathbb{E}^{R}\Big(\Gamma\int_0^t \langle \nabla f(\bar{X}_s), d\bar{K}_s\rangle\Big).
\end{eqnarray}
Indeed, the sequence $\Gamma\nabla f(\bar{X}^n_s)$ converges in distribution to $\Gamma\nabla f(\bar{X}_s)$ and the fact that $\lim_{a\uparrow+\infty}\sup_{n\geq1}R_n(\|K^n\|_{TV}>a)=0$ ensure that we can pass the limit in Eq.(\ref{equ5.12}). Combining (\ref{equ5.12}) with (\ref{equ5.11}), we show that $\bar{M}^{f,0}$ (which is the same as $\bar{M}^f$) is a solution of the martingale problem for Eq.(\ref{equ1.1}). Moreover, it is obvious that the condition $R\circ X_0^{-1}=\mathcal{L}_{X_0}$ hold.  $\Box$

\vskip0.3cm

\noindent{\bf Acknowledgements:}
Cheng L. is funded by the National Natural Science Foundation of China under Grants 12001272 and 12171173. Liu W. is supported by the National Natural Science Foundation of China under Grants 12471143 and 12131019.


\end{document}